\documentclass[11pt,a4paper]{article}
\usepackage[T1]{fontenc}
\usepackage{lmodern}
\usepackage{amsmath,amsthm,amssymb,mathtools}
\usepackage[margin=27mm,top=25mm,bottom=26mm,headheight=14pt]{geometry}
\usepackage{microtype,booktabs,xcolor,fancyhdr}
\definecolor{inkblue}{HTML}{253C54}
\usepackage[colorlinks=true,linkcolor=inkblue,citecolor=inkblue,urlcolor=inkblue]{hyperref}
\hypersetup{pdftitle={Sharp Gaussian approximation for Polyak--Ruppert averaged linear stochastic
approximation},pdfauthor={Bogdan Butyrin, Sergey Samsonov, Eric Moulines, Alexey Naumov},pdfsubject={Finite and limiting covariance targets}}
\newtheorem{theorem}{Theorem}[section]
\newtheorem{proposition}[theorem]{Proposition}
\newtheorem{lemma}[theorem]{Lemma}
\newtheorem{corollary}[theorem]{Corollary}
\theoremstyle{definition}

\newcommand{\dc}{\rho^{\mathrm{Conv}}}
\newcommand{\E}{\mathbb E}
\newcommand{\Pp}{\mathbb P}
\newcommand{\Cov}{\operatorname{Cov}}
\newcommand{\N}{\mathcal N}
\newcommand{\law}{\mathcal L}
\newcommand{\tr}{\operatorname{tr}}

\allowdisplaybreaks[2]
\title{Sharp Gaussian approximation for Polyak--Ruppert averaged linear
stochastic approximation}
\author{Bogdan Butyrin \and Eric Moulines \and Alexey Naumov \and Sergey Samsonov }
\date{11 November 2026}

\begin{document}
\maketitle

\begin{abstract}
We study how quickly the normalized Polyak--Ruppert average approaches
its Gaussian limit in linear stochastic approximation with bounded i.i.d.\
coefficients. For step sizes decaying as $k^{-\gamma}$, with
$1/2\le\gamma<1$, we prove an upper bound of order $n^{\gamma-1}$ on the
approximation error over convex sets in fixed dimension. A matching lower
bound holds for every fixed model. The fastest rate in this range is
$n^{-1/2}$, attained at $\gamma=1/2$. We also obtain an $O(n^{-1/2})$ bound on
the  error of the recursive multiplier bootstrap procedure.
\end{abstract}

\section{Introduction}\label{sec:intro}

Robbins and Monro \cite{rm} considered the problem of finding a root of an
unknown function from noisy evaluations. In vector form, let the mean
response at $\theta\in\mathbb R^d$ be
$h(\theta)=\mathbb E[F(\theta,Z_1)]$. The parameter of interest solves
\begin{equation}\label{eq:rm}
 h(\theta^\star)=0 .
\end{equation}
At iteration $k$, a fresh observation $Z_k$ provides the response
$F(\theta_{k-1},Z_k)$ at the current estimate. The Robbins--Monro update
subtracts $\alpha_k F(\theta_{k-1},Z_k)$ from $\theta_{k-1}$, where
$\alpha_k>0$ is the step size.

In linear stochastic approximation (LSA), the response takes the affine form
$F(\theta,z)=A(z)\theta-b(z)$. Writing $\bar A=\mathbb E[A(Z_1)]$ and
$\bar b=\mathbb E[b(Z_1)]$, equation~\eqref{eq:rm} reduces to
$\bar A\theta^\star=\bar b$. Starting from a deterministic $\theta_0$,
the recursion becomes
\begin{equation}\label{eq:lsa}
 \theta_k=\theta_{k-1}
     -\alpha_k\{A(Z_k)\theta_{k-1}-b(Z_k)\},\qquad k\ge1 .
\end{equation}
Ruppert \cite{ruppert} and Polyak \cite{polyak} proposed averaging the
iterates. For the average $\bar\theta_n=n^{-1}\sum_{k=0}^{n-1}\theta_k$,
the central limit theorem of Polyak and Juditsky \cite{pj} gives, under
appropriate stability, moment and step-size conditions,
\begin{equation}\label{eq:intro-clt}
 \sqrt n(\bar\theta_n-\theta^\star)
   \ \xrightarrow{\mathrm d}\ \mathcal N(0,\Sigma_\infty),\qquad
 \Sigma_\infty=\bar A^{-1}\Sigma_\varepsilon\bar A^{-\top}.
\end{equation}
Here $\Sigma_\varepsilon$ is the covariance of the response
$A(Z_1)\theta^\star-b(Z_1)$ at the root.

For $p\ge2$, the quantities $(\E\|\theta_k-\theta^\star\|^p)^{1/p}$
and $(\E\|\bar\theta_n-\theta^\star\|^p)^{1/p}$ measure the size of the
estimation error of the last and averaged iterates, respectively.
Durmus et al.\ \cite{durmus} derived bounds on
these moments and on deviation probabilities for constant step sizes,
under independent observations and under uniformly geometrically ergodic
Markov sampling. For polynomially decreasing step sizes, Butyrin
et al.\ \cite{butyrin} bounded the moments of the average, with leading
term proportional to $\sqrt{p\,\tr(\Sigma_\infty)/n}$.

A second line of work studies Gaussian approximation and statistical
inference. For confidence regions, the relevant error is the largest
difference between the sampling and Gaussian probabilities of a convex
Borel set. With independent observations, Samsonov et al.\ \cite{samsonov}
obtained a bound of order $n^{-1/4}$ up to logarithmic factors for
Polyak--Ruppert averaging. Wu et al.\ \cite{wu} obtained
order $n^{-1/3}$ for temporal difference learning with linear function
approximation. Butyrin et al.\ \cite{butyrin} proved this order for
general i.i.d. LSA. All three analyses use the normal approximation for
nonlinear statistics of independent observations due to Shao and Zhang
\cite{shao}.

To see where the exponent $1/3$ comes from, let $\Sigma_n$ denote the
covariance of the linear term in the expansion of
$\sqrt n(\bar\theta_n-\theta^\star)$, as defined in \eqref{eq:targets}.
For $\alpha_k=c_0(k+k_0)^{-\gamma}$, the bounds in \cite{wu,butyrin}
contain a term of order $n^{-\gamma/2}$ for approximation with this
finite-sample covariance and a term of order $n^{\gamma-1}$ for
replacing it by $\Sigma_\infty$. Balancing the two selects $\gamma=2/3$.
A scalar example with constant drift and additive Gaussian noise gives
a lower bound of order $n^{\gamma-1}$ at the limiting variance
\cite[eq.~(3.36)]{butyrin}.
For i.i.d. temporal difference learning, Cai et al.\ \cite[Theorem~3.9 and Corollary~3.10]{cai} obtained dimension-explicit
convex-distance bounds for the finite-sample Gaussian target. At fixed
dimension and for each fixed $1/2<\gamma<1$, their bounds (3.14) and
(3.15) are of order
$n^{-1/2}+n^{-\gamma/2}$, both for centering at the expectation of the
average and for centering at the root. For statistical inference in LSA,
Samsonov et al.\ \cite{samsonov} and Butyrin et al.\ \cite{butyrin}
obtained non-asymptotic bounds for the recursive multiplier bootstrap of Fang et al.\ \cite{fangboot}.

\paragraph{Main contributions.}
We improve the best exponent for the Polyak--Juditsky approximation of
general i.i.d. LSA from $1/3$ to $1/2$, attained at $\gamma=1/2$, and
give matching lower bounds. With the finite covariance $\Sigma_n$ the
error is $O(n^{-1/2})$ for every fixed $1/2\le\gamma<1$. A single
bounded scalar model attains this order even when the variance of the
centered Gaussian approximation is chosen to minimize the Kolmogorov distance.
At $\Sigma_\infty$ we prove a lower bound for every fixed model,
including at $\gamma=1/2$. For $1/2<\gamma<1$ we identify
the positive leading constant in the error of order $n^{\gamma-1}$.

The proof uses the expansion of the iterate error in successive
centered random matrix factors. The $L^2$ bound for the remainder after
the linear term contributes $O(n^{-\gamma/2})$ to the earlier
normal-approximation bounds. We retain several orders. Independence makes terms with
distinct sets of observation times orthogonal in the fluctuation started
at zero. After conditioning on the early observations, orthogonality
across distinct insertion times controls the error in truncating the
restarted recursion. We then apply the multivariate Stein coupling bound
of Fang and R\"ollin \cite{fr} to the retained terms as a function of the
remaining independent observations. For the bootstrap, comparing the sampling
law and the conditional bootstrap law to the same Gaussian law
$\N(0,\Sigma_n)$ gives an $O(n^{-1/2})$ bound on their expected convex
distance. The corresponding high-probability bound yields a coverage
guarantee for one-sided bootstrap intervals for $v^\top\theta^\star$.

\paragraph{Related work.}
For nonlinear SGD, Sheshukova et al.\ \cite{shesh} used the
Shao--Zhang approach \cite{shao} to prove bounds for Gaussian
approximation and for the multiplier bootstrap. Belomestny et al.\ \cite{bns} constructed a smooth, strongly convex scalar
example attaining the upper bound of order $n^{1/2-\gamma}+n^{\gamma-1}$
at the limiting variance for every $1/2<\gamma<1$. The term
$n^{1/2-\gamma}$ comes from curvature of the mean gradient and persists
at the finite linearized variance. For Markovian LSA, Samsonov et al.\ \cite{markovlsa} proved a bound of order
$n^{-1/4}$ up to logarithmic factors for scalar projections in
Kolmogorov distance and analyzed a multiplier block bootstrap. Wu et al.\ \cite{markovtd} obtained $O(n^{-1/4}\log n)$ over convex
sets for Markovian temporal difference learning. Under Markov sampling,
the current matrix can be correlated with the preceding iterate, so
the mean error recursion can contain a bias term. Under the i.i.d.
linear model \eqref{eq:lsa}, independence gives exactly
\[
 \E(\theta_k-\theta^\star)
   =(I_d-\alpha_k\bar A)\E(\theta_{k-1}-\theta^\star).
\]

\paragraph{Notation.}
We write $A_k=A(Z_k)$, $b_k=b(Z_k)$ and $\eta\sim\mathcal N(0,I_d)$.
The norm $\|\cdot\|$ is Euclidean for vectors and the induced operator
norm for matrices. For $Q=Q^\top\succ0$, let
$\|x\|_Q=(x^\top Qx)^{1/2}$ and
$\|B\|_Q=\sup_{x\ne0}\|Bx\|_Q/\|x\|_Q$.
We take $\prod_{\ell=m}^k B_\ell=B_k\cdots B_m$, with value $I_d$ for
an empty product. For probability laws $\mu,\nu$ on $\mathbb R^d$,
write
\[
 \dc(\mu,\nu)=
 \sup_{B\in\mathrm{Conv}(\mathbb R^d)}|\mu(B)-\nu(B)|,
\]
where $\mathrm{Conv}(\mathbb R^d)$ denotes the convex Borel sets.
Let $\Phi$ and $\varphi$ denote the standard normal distribution function and
density. For real-valued random variables $U,V$, the Kolmogorov distance is
$d_{\mathrm K}(U,V)=\sup_{t\in\mathbb R}|\Pp(U\le t)-\Pp(V\le t)|$.
We use the same distance notation for random variables and their laws.
For positive sequences $a_n,b_n$, write $a_n\lesssim b_n$ when
$a_n\le Cb_n$ for all sufficiently large $n$ with $C$ independent of $n$,
$a_n\asymp b_n$ when both $a_n\lesssim b_n$ and $b_n\lesssim a_n$,
and $a_n\sim b_n$ when $a_n/b_n\to1$.

\section{Main results}\label{sec:main}

We give rates in the sample size $n$. The dimension, the observation law, the
initial iterate and the step-size parameters are fixed. A statement for
all sufficiently large $n$ holds for $n\ge n_0$, where the threshold can
depend on these quantities. An explicit sufficient threshold can be
obtained from the inequalities in the proofs.

\subsection{Upper and lower bounds for Gaussian approximation}\label{sec:assumptions}
Let
$\widetilde A(z)=A(z)-\bar A$, $\widetilde b(z)=b(z)-\bar b$, and
$\varepsilon(z)=\widetilde A(z)\theta^\star-\widetilde b(z)$.
Thus $\varepsilon_k=\varepsilon(Z_k)=A_k\theta^\star-b_k$ and
$\widetilde A_k=\widetilde A(Z_k)$.

\paragraph{A1.} The observations $(Z_k)_{k\ge1}$ are i.i.d. on
$(\Omega,\mathcal F,\mathbb P)$ with common law $\pi$ on the measurable
space $(\mathsf Z,\mathcal Z)$.

\paragraph{A2.} The means of the measurable coefficient functions are
$\int A(z)\,\pi(dz)=\bar A$ and $\int b(z)\,\pi(dz)=\bar b$.
The matrix $-\bar A$ is Hurwitz. Moreover,
\[
 C_A:=\sup_{z\in\mathsf Z}\|A(z)\|\vee
        \sup_{z\in\mathsf Z}\|\widetilde A(z)\|<\infty,
 \qquad \|\varepsilon\|_\infty:=\sup_{z\in\mathsf Z}\|\varepsilon(z)\|<\infty,
\]
\[
 \Sigma_\varepsilon:=\int\varepsilon(z)\varepsilon(z)^\top\,\pi(dz),
 \qquad \lambda_{\min}:=\lambda_{\min}(\Sigma_\varepsilon)>0.
\]
In particular $\bar A$ is invertible, so the root $\theta^\star$ is unique.
For $P=P^\top\succ0$, let $Q=Q^\top\succ0$ be the solution of
$\bar A^\top Q+Q\bar A=P$. As in \cite[Proposition~1]{butyrin}, set
\begin{equation}\label{eq:lyapunovsetup}
 \kappa_Q=\frac{\lambda_{\max}(Q)}{\lambda_{\min}(Q)},\qquad
 \alpha_\infty=\frac{\lambda_{\min}(P)}{2\kappa_Q\|\bar A\|_Q^2}
          \wedge\frac{\|Q\|}{\lambda_{\min}(P)}.
\end{equation}

\paragraph{A3.} The step sizes are
\begin{equation}\label{eq:a3}
 \alpha_k=c_0(k+k_0)^{-\gamma},\qquad k\ge1,\qquad
 \tfrac12\le\gamma<1,\quad 0<c_0\le\alpha_\infty,\quad k_0\ge0,
\end{equation}
with $c_0$, $k_0$ and $\gamma$ fixed independently of $n$.

The offset $k_0$ controls the initial step sizes. For any fixed offset,
the steps eventually enter the range needed for the moment bounds, so
the finite initial segment affects only their constants.
The step-size condition in \cite[Assumption~A3$(p)$]{butyrin} uses
$1/2<\gamma<1$. To see why we can include $1/2$, observe that as $k\to\infty$,
\begin{equation}\label{eq:step-relative}
 \frac{\alpha_{k-1}/\alpha_k-1}{\alpha_k}
       \sim \frac{\gamma}{c_0}(k+k_0)^{\gamma-1}\longrightarrow0
       \qquad(1/2\le\gamma<1).
\end{equation}
Thus the relative change of the step size is smaller than the
one-step contraction of order $\alpha_k$. Lemma~\ref{lem:moments} uses
this comparison at $\gamma=1/2$, where $\sum_k\alpha_k^2$ diverges.
For the harmonic schedule $\gamma=1$, the ratio in \eqref{eq:step-relative}
tends to $1/c_0$ for every fixed $k_0$. The resulting balance between
step-size variation and contraction requires a separate analysis. 

Define the products and averaging weights by

\begin{align}
 G_{m:k}&=\prod_{\ell=m}^k(I_d-\alpha_\ell\bar A), &
 \Gamma_{m:k}&=\prod_{\ell=m}^k(I_d-\alpha_\ell A_\ell),\label{eq:prod-assumptions}\\
 Q_\ell&=\alpha_\ell\sum_{j=\ell}^{n-1}G_{\ell+1:j}, &1&\le\ell\le n-1 .\label{eq:weights-assumptions}
\end{align}
The horizon $n$ is suppressed in $Q_\ell$. The two covariance matrices are
\begin{equation}\label{eq:targets}
 \Sigma_n=\frac1n\sum_{\ell=1}^{n-1} Q_\ell\Sigma_\varepsilon Q_\ell^\top,
 \qquad
 \Sigma_\infty=\bar A^{-1}\Sigma_\varepsilon\bar A^{-\top}.
\end{equation}
Here $\Sigma_n$ is the covariance of the linear statistic
$-n^{-1/2}\sum_{\ell=1}^{n-1}Q_\ell\varepsilon_\ell$. The effect of
the higher-order terms on the covariance of the full average is quantified below.

\paragraph{A4.} The sample size is large enough that
\begin{equation}\label{eq:a4}
 \lambda_{\min}(\Sigma_n)\ge\tfrac12\lambda_{\min}(\Sigma_\infty).
\end{equation}
Since $\Sigma_\infty\succ0$, the convergence in Lemma~\ref{lem:matrix}
ensures A4 for all sufficiently large $n$.

Write $T_n=\sqrt n(\bar\theta_n-\theta^\star)$.

\begin{theorem}\label{thm:finite}
Under A1--A4 there is a constant $C<\infty$ such that, for all
sufficiently large $n$,
\begin{align}
 \dc(T_n,\Sigma_n^{1/2}\eta)&\le C n^{-1/2},\label{eq:finiteconv}\\
 \dc(T_n,\Sigma_\infty^{1/2}\eta)&\le C n^{\gamma-1}.
 \label{eq:infinity-upper}
\end{align}
\end{theorem}
The constant $C$ is independent of $n$. Its dependence on the dimension
and on the model and schedule parameters can be tracked through the
moment estimates, the covariance eigenvalue bounds and the multivariate
Stein inequality. The explicit dimension factors in the latter are
powers of $d$ given in \eqref{eq:fr-bound}. We suppress these constants
to state the rates in $n$.

Let $L$ be the unique symmetric solution of the Lyapunov equation
\begin{equation}\label{eq:lyap}
 c_0(\bar A L+L\bar A^\top)=\Sigma_\infty,\qquad
 B=\Sigma_\infty^{-1/2}L\Sigma_\infty^{-1/2}.
\end{equation}
Since $-\bar A$ is Hurwitz and $\Sigma_\infty\succ0$, the integral
$L=c_0^{-1}\int_0^\infty e^{-t\bar A}\Sigma_\infty
e^{-t\bar A^\top}\,dt$ is positive definite. For the averaging weights in
\eqref{eq:weights-assumptions}, Lemma~\ref{lem:matrix} gives the expansion
in operator norm
\begin{equation}\label{eq:exp}
 \Sigma_n=\Sigma_\infty+n^{\gamma-1}L+o(n^{\gamma-1}),
 \qquad \tfrac12\le\gamma<1 .
\end{equation}
Wu et al.\ \cite[Theorem~6.11]{wu} obtained the corresponding expansion
for temporal difference learning.
\begin{theorem}
\label{thm:limit}
Under A1--A4, for every fixed $1/2\le\gamma<1$,
\[
 \lim_{n\to\infty} n^{1-\gamma}
      \dc(\Sigma_n^{1/2}\eta,\Sigma_\infty^{1/2}\eta)
       = C_*,\qquad C_*:=\tfrac14\E|\eta^\top B\eta-\tr B|>0.
\]
Furthemore,
\begin{equation}\label{eq:limit-order}
 \dc(T_n,\Sigma_\infty^{1/2}\eta)\asymp n^{\gamma-1},
 \qquad \tfrac12\le\gamma<1 .
\end{equation}
If $1/2<\gamma<1$, its leading constant is also $C_*$:
\begin{equation}\label{eq:exactrate}
 \lim_{n\to\infty} n^{1-\gamma}
       \dc(T_n,\Sigma_\infty^{1/2}\eta)=C_* .
\end{equation}
\end{theorem}
At $\gamma=1/2$, \eqref{eq:limit-order} gives the lower bound
$c n^{-1/2}$ for each fixed model. The constant in
\eqref{eq:exactrate} need not apply at this endpoint, as the
initialization example at the end of the proof shows.

The next proposition establishes sharpness of the finite-sample bound
\eqref{eq:finiteconv} over the model class.

\begin{proposition}\label{prop:sharp-example}
Let $A_k=1+R_k/2$, $b_k=\xi_k$, where the $R_k$ are independent
Rademacher variables, independent of the entire i.i.d. sequence $(\xi_k)$,
and
\[
 \Pp(\xi_k=\sqrt{2})=1/3,\qquad \Pp(\xi_k=-1/\sqrt{2})=2/3 .
\]
Take $\theta_0=\theta^\star=0$, $c_0=1/4$, and fix $k_0\ge0$ and
$1/2\le\gamma<1$. There are constants $0<c\le C<\infty$ such that,
for all sufficiently large $n$,
\begin{equation}\label{eq:examplefinite}
 \frac{c}{\sqrt n}
 \le \inf_{\sigma>0}d_{\mathrm K}(T_n,\sigma\eta)
 \le d_{\mathrm K}(T_n,\Sigma_n^{1/2}\eta)
 \le \frac{C}{\sqrt n}.
\end{equation}
\end{proposition}
In one dimension, half-lines are convex and every interval probability
is a difference of two distribution function values or their one-sided
limits. Hence $d_{\mathrm K}(U,V)\le\dc(U,V)\le 2d_{\mathrm K}(U,V)$.
The lower bound in \eqref{eq:examplefinite} therefore also applies to
the convex distance in \eqref{eq:finiteconv}.

\section{Recursive multiplier bootstrap}\label{sec:applications}
For each replicate $\ell$, take weights $(w_k^{\ell})_{k\ge1}$ independent across times
and replicates and independent of the data. With $\theta_0^{b,\ell}=\theta_0$, compute
\begin{equation}\label{eq:bootstrap-procedure}
 \theta_k^{b,\ell}=\theta_{k-1}^{b,\ell}
  -\alpha_k w_k^{\ell}(A_k\theta_{k-1}^{b,\ell}-b_k),
 \qquad \bar\theta_n^{b,\ell}=\frac1n\sum_{k=0}^{n-1}\theta_k^{b,\ell}.
\end{equation}
This is the online recursion of \cite{fangboot,samsonov,butyrin}.
Write $Z^{n-1}=(Z_1,\ldots,Z_{n-1})$ and let $\mathbb P^b,\mathbb E^b$ denote
probability and expectation conditional on $Z^{n-1}$. Dropping the replicate index,
set $U_n^b=\sqrt n(\bar\theta_n^b-\bar\theta_n)$ and
\begin{equation}\label{eq:bootstrap-distance}
 D_n=\sup_{B\in\mathrm{Conv}(\mathbb R^d)}
 \left|\mathbb P^b(U_n^b\in B)-\mathbb P\{\sqrt n(\bar\theta_n-\theta^\star)\in B\}\right|.
\end{equation}
The reference law in \eqref{eq:bootstrap-distance} is the unconditional sampling
law. The bounds below concern the exact conditional law of one replicate.

\begin{theorem}
\label{thm:boot}
Assume A1--A4 and fix $d$. Suppose the weights are i.i.d., almost surely bounded, with
$\mathbb Ew_1=1$ and $\operatorname{Var}(w_1)=1$. For all sufficiently large $n$,
\[
 \mathbb E D_n\le C n^{-1/2},\qquad
 \mathbb P\!\left(D_n\le C\sqrt{\frac{\log n}{n}}\right)\ge1-\frac1n .
\]
The constant may depend on the model, schedule, dimension and multiplier bound.
\end{theorem}

Weights taking the values $0$ and $2$ with equal probability satisfy the theorem.
In the proof the conditional covariance is compared directly with
$\Sigma_n$, which avoids the $n^{\gamma-1}$ error of replacing it by $\Sigma_\infty$.

\begin{corollary}
\label{cor:bootcoverage}
Under Theorem~\ref{thm:boot}, fix $v\ne0$ and $0<\delta<1$. Let
$q^b_{1-\delta}(v)=\inf\{t:\mathbb P^b(v^\top U_n^b\le t)\ge1-\delta\}$.
Then
\[
 \left|\mathbb P\!\left\{
 v^\top\theta^\star\ge v^\top\bar\theta_n-
             \frac{q^b_{1-\delta}(v)}{\sqrt n}\right\}-(1-\delta)\right|
 \le C\sqrt{\frac{\log n}{n}} .
\]
\end{corollary}

\section{Proofs of the main results}\label{sec:proofs}

\subsection{Proof of Theorem~\ref{thm:finite}}\label{sec:proof-finite}

\paragraph{Conditioning on the first $m$ observations.}
Write $e_k=\theta_k-\theta^\star$, $\mathcal F_k=\sigma(Z_1,\ldots,Z_k)$,
and $\mathcal F_0=\{\varnothing,\Omega\}$. With
$Y_j=\varepsilon_j+\widetilde A_j e_{j-1}$, independence gives
$\E[Y_j\mid\mathcal F_{j-1}]=0$. Iterating the recursion around $\bar A$
and interchanging the two time sums gives
\begin{equation}\label{eq:mart}
 T_n=a_n-\frac1{\sqrt n}\sum_{j=1}^{n-1}Q_jY_j,\qquad
 a_n=\E T_n=\frac1{\sqrt n}\sum_{k=0}^{n-1}G_{1:k}e_0 .
\end{equation}
Choose $m=\lfloor\sqrt n\rfloor$. After time $m$ the step sizes are at
most $\alpha_m$, which makes successive terms of the expansion small.
The early observations enter through the starting value $e_m$ and
the conditional mean $\mu_m=\E[T_n\mid\mathcal F_m]$. Their contribution
after centering is controlled by
\begin{equation}\label{eq:prefix-variance}
 \E\|\mu_m-a_n\|^2
 =\frac1n\sum_{j=1}^m\E\|Q_jY_j\|^2\le C\frac mn .
\end{equation}
Indeed, conditioning \eqref{eq:mart} retains just the increments
through $m$, and distinct martingale increments are orthogonal.
Lemmas~\ref{lem:moments} and~\ref{lem:kernel} bound their second
moments and give $\|a_n\|\le Cn^{-1/2}$ by summability of $\|G_{1:k}\|$.

For deterministic $x\in\mathbb R^d$, denote by $e_k^{m,x}$ the
recursion started from $e_m^{m,x}=x$ and driven by
$X=(Z_{m+1},\ldots,Z_{n-1})$. Under $\Pp_x$ the $N=n-1-m$ inputs in $X$
have law $\pi^{\otimes N}$, and $\E_x$ denotes integration over these
inputs. Set
\[
 P_m=\frac1{\sqrt n}\sum_{k=0}^{m-1}e_k,\qquad
 F_x=\frac1{\sqrt n}\sum_{k=m}^{n-1}e_k^{m,x},\qquad Z_x=F_x-\E_x F_x .
\]
Then $T_n=P_m+F_{e_m}$. By independence, $\E_x e_k^{m,x}=G_{m+1:k}x$.
Thus
\begin{equation}\label{eq:sketch-condlaw}
 \mu_m=P_m+\frac1{\sqrt n}\sum_{k=m}^{n-1}G_{m+1:k}e_m,\qquad
 \law(T_n-\mu_m\mid\mathcal F_m)
   =\left.\law_{\Pp_x}(Z_x)\right|_{x=e_m}.
\end{equation}

\paragraph{Reduction to a deterministic starting value.}
Let $\nu_n=\N(0,\Sigma_n)$ and put
$d_{m,n}=\sup_{\|x\|\le1}\dc(\law_{\Pp_x}(Z_x),\nu_n)$.
Fix a convex Borel set $D$ and write $g_D(u)=\nu_n(D-u)$.
For a realization of the first $m$ observations, both $e_m$ and $\mu_m$
are fixed, so \eqref{eq:sketch-condlaw} gives
\[
 \Pp(T_n\in D\mid\mathcal F_m)
 =\left.\Pp_x(Z_x\in D-u)\right|_{(x,u)=(e_m,\mu_m)} .
\]
On $\{\|e_m\|\le1\}$ the difference from $g_D(\mu_m)$ is at most
$d_{m,n}$ because $D-\mu_m$ is convex for each realized translation.
On the complementary event the difference is bounded by one. Taking
expectations for this fixed $D$ yields
\[
 |\Pp(T_n\in D)-\E g_D(\mu_m)|
 \le d_{m,n}+\Pp(\|e_m\|>1)
 \le d_{m,n}+C\alpha_m^2 .
\]
The last inequality uses the fourth-moment estimate in Lemma~\ref{lem:moments}.
To integrate the translation itself, apply Taylor's formula at $a_n$.
The linear term vanishes since $\E(\mu_m-a_n)=0$. A4, the uniform upper
bound on $\|\Sigma_n\|$ from Lemma~\ref{lem:kernel}, and the Gaussian
derivative bounds in Lemma~\ref{lem:gausstools} yield
\[
 \begin{aligned}
 |\E g_D(\mu_m)-g_D(a_n)|
    &\le \tfrac12\sup_u\|\nabla^2g_D(u)\|\,\E\|\mu_m-a_n\|^2
      \le Cm/n,\\
 |g_D(a_n)-\nu_n(D)|
    &\le \sup_u\|\nabla g_D(u)\|\,\|a_n\|\le Cn^{-1/2}.
 \end{aligned}
\]
All constants are independent of $D$. We can now take the supremum over $D$.
Since $m/n\le n^{-1/2}$ and $\alpha_m^2=O(n^{-\gamma})$, we obtain
\begin{equation}\label{eq:prefix-reduction}
 \dc(T_n,\nu_n)\le d_{m,n}+Cn^{-1/2}.
\end{equation}

\paragraph{The expansion at the restart.}
For a fixed $x$ write $J_k^{(\ell)}=J_{m,k}^{(\ell)}(x)$ and
$H_k^{(K)}=H_{m,k}^{(K)}(x)$. Start with $J_m^{(0)}=x$ and
$J_m^{(\ell)}=0$ for $\ell\ge1$. For $k>m$ put
\begin{align}
 J_k^{(0)}&=(I_d-\alpha_k\bar A)J_{k-1}^{(0)}-\alpha_k\varepsilon_k,\notag\\
 J_k^{(\ell)}&=(I_d-\alpha_k\bar A)J_{k-1}^{(\ell)}
      -\alpha_k\widetilde A_k J_{k-1}^{(\ell-1)},\qquad\ell\ge1.
 \label{eq:J}
\end{align}
The remainder $H_k^{(K)}=e_k^{m,x}-\sum_{\ell=0}^KJ_k^{(\ell)}$
starts from zero and satisfies
\begin{equation}\label{eq:H}
 H_k^{(K)}=(I_d-\alpha_k A_k)H_{k-1}^{(K)}
                    -\alpha_k\widetilde A_k J_{k-1}^{(K)} .
\end{equation}
In particular, taking $m=0$ and $x=0$ recovers the fluctuation expansion
of Butyrin et al.\ \cite{butyrin}. Adding the
homogeneous solution for the original deterministic $e_0$ gives
\begin{equation}\label{eq:global-decomp}
 e_k=\Gamma_{1:k}e_0+\sum_{\ell=0}^K J_{0,k}^{(\ell)}(0)+H_{0,k}^{(K)}(0).
\end{equation}
At zero initialization the normalized sum of the zeroth-order terms
is $W_n=-n^{-1/2}\sum_{j<n}Q_j\varepsilon_j$. For this zero-start trajectory, \eqref{eq:mart} and orthogonality give
\[
 \E\|T_n-W_n\|^2
 =\frac1n\sum_{j=1}^{n-1}\E\|Q_j\widetilde A_j e_{j-1}\|^2
 \lesssim n^{-\gamma}.
\]
In the zero-start expansion each order $\ell$ contains products indexed
by $\ell+1$ distinct observation times. Two different time sets have a
time occurring in just one product. Integrating over the observation at
that time gives zero, since its factor is centered. After a nonzero
restart we use orthogonality over the time of the last insertion, as follows.

For any $m$ and deterministic $x$ define
\begin{equation}\label{eq:sketch-FH}
 F_x^{[K]}=\frac1{\sqrt n}\sum_{k=m}^{n-1}\sum_{\ell=0}^{K}J_{m,k}^{(\ell)}(x),
 \qquad F_x-F_x^{[K]}=\frac1{\sqrt n}\sum_{k=m}^{n-1}H_{m,k}^{(K)}(x).
\end{equation}
We now take $m=\lfloor\sqrt n\rfloor$ and $\|x\|\le1$.
Each additional order inserts one centered matrix. More precisely,
for an adapted sequence $v$ let
$(\mathcal T_m v)_k=-\sum_{j=m+1}^{k}\alpha_jG_{j+1:k}\widetilde A_jv_{j-1}$.
In component form, \eqref{eq:operatorbridge} gives
\[
 J_{m,k}^{(\ell+1)}(x)=(\mathcal T_m J_{m,\bullet}^{(\ell)}(x))_k,\qquad
 H_{m,k}^{(K)}(x)=(\mathcal T_m^{K+1}e_\bullet^{m,x})_k ,
\]
where $\bullet$ denotes the sequence of time indices from $m$ to $n-1$. At fixed $k$, orthogonality over the insertion time $j$ gives
\[
 \E_x\|(\mathcal T_m v)_k\|^2
   =\sum_{j=m+1}^k\alpha_j^2\E_x\|G_{j+1:k}\widetilde A_jv_{j-1}\|^2
   \le C\alpha_m\max_{m\le t<n}\E_x\|v_t\|^2 .
\]
Thus Lemma~\ref{lem:proxy} gives, for $K=15$ and $F=F_x^{[15]}$,
\[
 \E_xF=\E_xF_x,\qquad
 \|F_x-F\|_{L^2(\Pp_x)}\le C\sqrt n\,\alpha_m^8
 \le C n^{1/2-4\gamma}\le Cn^{-3/2}.
\]

\paragraph{Normal approximation under $\Pp_x$.}
Write $V_x=\Cov_x(F_x)$, $V_x^{[15]}=\Cov_x(F)$ and
$Z_x^{[15]}=F-\E_x F$. Lemma~\ref{lem:suffixcov} gives
\[
 \|V_x-\Sigma_n\|\le C\{m/n+n^{-\gamma}+(n\alpha_m)^{-1}\}
                     \le Cn^{-1/2}.
\]
Hence $\E_x\|Z_x\|^2=\tr V_x\le C$. The truncation bound and
Cauchy--Schwarz give
\[
 \|V_x^{[15]}-V_x\|
 \le (\|Z_x^{[15]}\|_{L^2(\Pp_x)}+\|Z_x\|_{L^2(\Pp_x)})
             \|Z_x^{[15]}-Z_x\|_{L^2(\Pp_x)}
 \le Cn^{-3/2}.
\]
By A4 it follows that $\lambda_{\min}(V_x^{[15]})\ge
\lambda_{\min}(\Sigma_\infty)/4$ for all sufficiently large $n$,
uniformly in $\|x\|\le1$. For the Fang--R\"ollin bound, we normalize
the centered sum $Z_x^{[15]}$ by multiplying it by $(V_x^{[15]})^{-1/2}$.
The resulting vector has covariance $I_d$ under $\Pp_x$.

It remains to check the independent-coordinate criterion in
Lemma~\ref{lem:sensitivity} under the product law of $X$.
Use the replacements in \eqref{eq:resample-differences} for coordinates
$i,j\in\{m+1,\ldots,n-1\}$ and let $s_k=\sum_{t=1}^k\alpha_t$. A change at $i$ enters the retained
iterates with a factor $\alpha_i$. Its contribution at time $k\ge i$
is bounded by $C\alpha_i e^{-b(s_k-s_i)}$. A further change at $j>i$
acts on the perturbation reaching $j-1$. Its mixed contribution at time
$k\ge j$ is bounded by
$C\alpha_i\alpha_j e^{-b(s_{j-1}-s_i)}e^{-b(s_k-s_j)}$.
Summation over $k$ costs at most $C/\alpha_i$ and $C/\alpha_j$,
respectively, by Lemma~\ref{lem:kernel}. Division by $\sqrt n$ therefore yields
\begin{equation}\label{eq:sketch-influences}
 \begin{aligned}
 \|\Delta_i F\|&\le Cn^{-1/2},&
 \|\Delta_i\Delta_j F\|&\le Cn^{-1/2}h_{ij}\quad(i\ne j),\\
 h_{ii}&=1,& h_{ij}=h_{ji}&=\alpha_i e^{-b(s_{j-1}-s_i)}\quad(i<j),
 \end{aligned}
 \qquad \sup_i\sum_j h_{ij}\le C .
\end{equation}
Centering $F$ preserves these differences. Since $N=n-1-m\asymp n$,
the bounds imply the $N^{-1/2}$ hypotheses of Lemma~\ref{lem:sensitivity}.
That lemma and Gaussian comparison in Lemma~\ref{lem:gausstools} give
\[
 \dc(\law_{\Pp_x}(Z_x^{[15]}),\nu_n)
 \le C N^{-1/2}+C\|V_x^{[15]}-\Sigma_n\|\le Cn^{-1/2}.
\]

\paragraph{Restoring the exact remainder.}
Put $t=n^{-1/2}$. Since $Z_x-Z_x^{[15]}=F_x-F$, Markov's inequality
bounds $p_x:=\Pp_x(\|Z_x-Z_x^{[15]}\|>t)$ by
$t^{-2}\E_x\|F_x-F\|^2\le Cn^{-2}$. For any convex Borel set $D$,
its inner and outer $t$-neighborhoods satisfy
\[
 \Pp_x(Z_x^{[15]}\in D^{-t})-p_x
 \le \Pp_x(Z_x\in D)
 \le \Pp_x(Z_x^{[15]}\in D^t)+p_x .
\]
Applying the preceding normal approximation to each neighborhood and
the Gaussian boundary bound of Lemma~\ref{lem:gausstools} gives
\begin{equation}\label{eq:restart-centered-normal}
 |\Pp_x(Z_x\in D)-\nu_n(D)|\le Cn^{-1/2}+Ct+p_x\le C'n^{-1/2}.
\end{equation}
The bound is uniform in $D$ and $\|x\|\le1$, so it gives
$d_{m,n}\le C'n^{-1/2}$. Equation~\eqref{eq:prefix-reduction} proves
\eqref{eq:finiteconv}. Finally, Lemmas~\ref{lem:gausstools}
and~\ref{lem:matrix} give
\[
 \dc(T_n,\Sigma_\infty^{1/2}\eta)
   \le \dc(T_n,\nu_n)+C\|\Sigma_n-\Sigma_\infty\|
   \le C(n^{-1/2}+n^{\gamma-1}).
 \tag*{\(\square\)}
\]

\subsection{Proof of Theorem~\ref{thm:limit}}

\paragraph{The Gaussian probabilities.}
Write $h_n=n^{\gamma-1}$. To normalize the limiting covariance to $I_d$,
apply $z\mapsto\Sigma_\infty^{-1/2}z$ to both Gaussian laws.
This invertible map preserves convex distance. After this normalization,
the finite-sample Gaussian law has covariance
\[
 V_n=\Sigma_\infty^{-1/2}\Sigma_n\Sigma_\infty^{-1/2}=I_d+h_n B+o(h_n)
\]
by \eqref{eq:exp}. Since $B\succ0$, $V_n-I_d\succ0$ for large $n$.
Let $p_V$ denote the $\N(0,V)$ density. The sign of $p_{I_d}-p_{V_n}$
changes on the boundary of the ellipsoid
\[
 D_n=\{z:z^\top(I_d-V_n^{-1})z\le\log\det V_n\}.
\]
Indeed this inequality is exactly $p_{I_d}(z)\ge p_{V_n}(z)$. Integrating
the positive part over $D_n$ gives half the $L^1$ distance between the
densities, because the densities have equal integrals. No Borel set can
give a larger absolute probability difference, and $D_n$ is convex.
Thus Lemma~\ref{lem:ordered-gauss} gives
\[
 \dc(\Sigma_n^{1/2}\eta,\Sigma_\infty^{1/2}\eta)
   =\tfrac12\|p_{V_n}-p_{I_d}\|_1 .
\]

\paragraph{The first-order coefficient.}
Differentiating the Gaussian density at $I_d$ in direction $B$ gives
\[
 \left.\frac{d}{dt}p_{I_d+tB}(z)\right|_{t=0}
   =\tfrac12 p_{I_d}(z)(z^\top Bz-\tr B).
\]
For $V$ near $I_d$, the first and second derivatives in matrix directions
of norm at most one are bounded by $C(1+\|z\|^4)e^{-c\|z\|^2}$.
Taylor's formula and $V_n-I_d=h_nB+o(h_n)$ therefore give $L^1$ convergence:
\[
 \frac{p_{V_n}-p_{I_d}}{h_n}
 \ \longrightarrow\ \tfrac12p_{I_d}(z)(z^\top Bz-\tr B).
\]
Taking norms proves the Gaussian equivalent, with
$C_*=\tfrac14\E|\eta^\top B\eta-\tr B|>0$ since
$\operatorname{Var}(\eta^\top B\eta)=2\tr(B^2)>0$.

\paragraph{Transfer to the averaged iterate.}
By the reverse triangle inequality and Theorem~\ref{thm:finite},
\[
 \left|\dc(T_n,\Sigma_\infty^{1/2}\eta)
       -\dc(\Sigma_n^{1/2}\eta,\Sigma_\infty^{1/2}\eta)\right|
 \le \dc(T_n,\Sigma_n^{1/2}\eta)\le Cn^{-1/2}.
\]
For $1/2<\gamma<1$ we have $n^{-1/2}=o(h_n)$, which proves
\eqref{eq:exactrate}. For a fixed contrast $v\ne0$, set
$s_v^2=v^\top\Sigma_\infty v$. Applying the same bound to the halfspace
$\{z:v^\top z\le s_v\}$ and expanding its Gaussian probability gives
\[
 \begin{aligned}
 \Pp(v^\top T_n\le s_v)
 &=\Phi\!\left(\frac{s_v}{\sqrt{v^\top\Sigma_n v}}\right)+O(n^{-1/2})\\
 &=\Phi(1)-h_n\frac{\varphi(1)v^\top Lv}{2s_v^2}+o(h_n).
 \end{aligned}
\]

\paragraph{The lower bound at $\gamma=1/2$.}
Fix a unit vector $v$ and put $s_v^2=v^\top\Sigma_\infty v$,
$\ell_v=v^\top Lv>0$, $m_n=v^\top a_n$ and
$\tau_n^2=\operatorname{Var}(v^\top T_n)$. In \eqref{eq:covexact} the
terms $\mathcal V(M_{j-1})$ are positive semidefinite, while
$\sup_j\|Q_j\|\le C$ and $\sum_{j<n}\|\mu_{j-1}\|\le C$.
Thus \eqref{eq:exp} and \eqref{eq:covexact} imply, for large $n$,
\[
 \tau_n^2\ge v^\top\Sigma_n v-C/n
       \ge s_v^2+\frac{\ell_v}{2\sqrt n}.
\]
To detect this excess use $f_R(x)=\exp(-x^2/(2R^2))$, $R\ge1$.
For $Z\sim\N(0,1)$ and $u\ge0$,
\[
 \E f_R(m+\sqrt u Z)=g_R(u)\exp\!\left(-\frac{m^2}{2(R^2+u)}\right),
 \qquad g_R(u)=(1+u/R^2)^{-1/2}.
\]
On the bounded range of $s_v^2$ and $\tau_n^2$, $-g_R'(u)\ge cR^{-2}$
uniformly for $R\ge1$. Since $\|f_R'''\|_\infty\le C R^{-3}$,
Lemma~\ref{lem:endpoint-smooth} gives
\[
 \E f_R(s_vZ)-\E f_R(v^\top T_n)
 \ge g_R(s_v^2)-g_R(\tau_n^2)-\frac{C}{R^3\sqrt n}
 \ge \frac{1}{\sqrt n}\left(\frac{c\ell_v}{R^2}-\frac{C}{R^3}\right).
\]
Fix $R$ large enough to make the bracket positive. For $0<u<1$, the
set $\{z:f_R(v^\top z)>u\}$ is a convex slab. Integrating its
probability over $u\in(0,1)$ bounds the left side by
$\dc(T_n,\Sigma_\infty^{1/2}\eta)$. The upper bound follows from
\eqref{eq:infinity-upper}, proving \eqref{eq:limit-order}.

\hfill$\square$


\subsection{Proof of Proposition~\ref{prop:sharp-example}}

\paragraph{Conditioning on the coefficients.}
In this scalar model $\bar A=\Sigma_\varepsilon=1$ and
$\varepsilon_i=-\xi_i$. With $P=1$ we have $Q=\alpha_\infty=1/2$,
so $c_0=1/4$ satisfies A3. Since $\theta_0=0$, the exact recursion gives
\[
 T_n=\frac1{\sqrt n}\sum_{i=1}^{n-1}U_i\xi_i,\qquad
 U_i=\alpha_i\sum_{k=i}^{n-1}\prod_{j=i+1}^k(1-\alpha_j A_j).
\]
Conditionally on $\mathcal A=\sigma(A_1,\ldots,A_{n-1})$ the weights are
fixed, and the $\xi_i$ retain their independent common law. Each factor
is positive since $\alpha_jA_j\le3/8$. The product is at most
$\exp\{-(s_k-s_i)/2\}$, with $s_k=\sum_{j=1}^k\alpha_j$, so
Lemma~\ref{lem:kernel} gives $U_i\le C$. For $n/4\le i\le n/2$,
the first $q_i=\lfloor1/\alpha_i\rfloor$ summands are within the horizon
for large $n$, since $q_i=O(n^\gamma)=o(n)$. Consequently
\[
 U_i\ge\alpha_iq_i(1-3\alpha_i/2)^{q_i}\ge c>0 .
\]
It follows that $v_n=n^{-1}\sum_{i<n}U_i^2$ and
$m_{3,n}=n^{-1}\sum_{i<n}U_i^3$ lie in a fixed compact subinterval of
$(0,\infty)$, uniformly in the realized coefficients.

\paragraph{A lower bound against any centered Gaussian.}
Here $\E\xi_i=0$, $\E\xi_i^2=1$ and $\kappa_3=\E\xi_i^3=1/\sqrt2$.
For $\chi(u)=\E e^{iu\xi_i}$, Taylor expansion around zero gives
$\log\chi(u)=-u^2/2-i\kappa_3u^3/6+O(u^4)$, with the logarithm chosen
continuously at zero. Conditional independence now yields, uniformly in
$1\le t\le2$ and in $\mathcal A$,
\[
 \begin{aligned}
 \E[e^{itT_n}\mid\mathcal A]
   &=\prod_{i<n}\chi(tU_i/\sqrt n)\\
   &=e^{-t^2v_n/2}
      \left(1-\frac{i\kappa_3t^3m_{3,n}}{6\sqrt n}+O(n^{-1})\right).
 \end{aligned}
\]
Indeed $\sum_{i<n}U_i^4/n^2=O(n^{-1})$. Thus
$\E[\sin(tT_n)\mid\mathcal A]\le-c/\sqrt n$ on this entire interval
for large $n$. Take a nonnegative $\psi\in C_c^\infty((1,2))$ with
$\int\psi=1$ and put $h(x)=\int\psi(t)\sin(tx)\,dt$. The function is
odd, and $h'$ is integrable by integration by parts in $t$. Taking expectations over the
coefficients gives $\E h(T_n)\le -c/\sqrt n$. For every $\sigma>0$,
symmetry gives $\E h(\sigma\eta)=0$. Integration by parts against the two
distribution functions yields
\[
 \frac{c}{\sqrt n}
 \le |\E h(T_n)-\E h(\sigma\eta)|
 \le \|h'\|_1\,d_{\mathrm K}(T_n,\sigma\eta).
\]
Taking the infimum over $\sigma>0$ proves the first inequality in
\eqref{eq:examplefinite}. By A4, $\Sigma_n>0$ for large $n$, so
$\sigma=\sqrt{\Sigma_n}$ is an admissible choice. Finally,
$d_{\mathrm K}\le\dc$ and \eqref{eq:finiteconv} give the upper bound.
\hfill$\square$

\section*{Declaration of AI use}
The authors used ChatGPT (OpenAI) for checking arguments for correctness, and for polishing the text of the manuscript. The research questions and proposed method were developed by the authors. The authors take full responsibility for the content of the paper, including the correctness of all results, proofs, and references.

\clearpage\appendix

\section{Auxiliary lemmas}\label{app:inputs}

\subsection{Moment bounds for the LSA recursion}
Let $s_0=0$, $s_k=\sum_{j=1}^k\alpha_j$ and
$\|Y\|_{L^p}=(\E\|Y\|^p)^{1/p}$, with the Euclidean norm for a vector
and the induced operator norm for a matrix. The products $G_{i:k}$ and
$\Gamma_{i:k}$ are defined in \eqref{eq:prod-assumptions}.

\begin{lemma}
\label{lem:moments}
Assume A1--A3, let $\theta_0$ be deterministic and fix $p\ge2$.
There are $C_p<\infty$ and $b_p>0$, depending on the law, $d$, $p$
and the fixed schedule, such that
\begin{align}
 \|G_{i:k}\|+\|\Gamma_{i:k}\|_{L^p}
    &\le C_p e^{-b_p(s_k-s_{i-1})} &&(1\le i\le k),\notag\\
 \|\theta_k-\theta^\star\|_{L^p}
    &\le C_p\{e^{-b_ps_k}\|\theta_0-\theta^\star\|+
                 \sqrt{\alpha_k}\} &&(k\ge1),\notag\\
 \|T_n\|_{L^p}
    &\le C_p(1+\|\theta_0-\theta^\star\|) &&(n\ge1).
 \label{eq:inputs}
\end{align}
\end{lemma}
For $1/2<\gamma<1$, the product, iterate and average estimates are given in
\cite[Proposition~1, Lemma~1, Proposition~2 and Theorem~1]{butyrin}.
The proofs for the random product and last iterate are in Appendix~B
of that paper, and the averaged bound is proved in Section~5.1.
The last-iterate estimate also follows from the argument of
\cite[Proposition~4]{samsonov}.
To use these fixed-order estimates with A3, start at a deterministic
time $j_0$ sufficiently large for the chosen $p$. The subsequent steps
form the same polynomial schedule with shifted time origin. Conditional
on $\theta_{j_0}$, the subsequent observations are independent with
law $\pi$. Bounded inputs give
$\|\theta_{j_0}-\theta^\star\|\le C_{j_0}(1+\|\theta_0-\theta^\star\|)$
almost surely.
Equivalently one may apply the fluctuation estimate at zero initial
error and propagate the deterministic initial error with $\Gamma_{1:k}$.
Including the finite initial products and the initial sum in the average
only enlarges $C_p$. The same argument applies to bounded weighted
i.i.d. inputs with the same mean matrix.

\begin{proof}[Moment bounds at $\gamma=1/2$]
First take an even integer $p\ge2$.
\smallskip\noindent\emph{Obtain a one-step contraction.}
Let $e_k=\theta_k-\theta^\star$, $\alpha=\alpha_k$ and
$y=-\alpha(A_kx+\varepsilon_k)$ when $e_{k-1}=x$. With
$\bar A^\top Q+Q\bar A=P\succ0$,
\[
 \E\!\left[p\|x\|_Q^{p-2}\,x^\top Qy\mid e_{k-1}=x\right]
 =-\tfrac p2\alpha\|x\|_Q^{p-2}x^\top P x
 \le -b'_p\alpha\|x\|_Q^p.
\]
At $x=0$ the derivative is zero. The Taylor remainder for
$\|x+y\|_Q^p$ is bounded by
$C_p(\|x\|_Q^{p-2}\|y\|_Q^2+\|y\|_Q^p)$.
Since $\|y\|_Q\le C\alpha(1+\|x\|_Q)$, for all sufficiently late $k$
this yields
\[
 \E[\|e_k\|_Q^p\mid e_{k-1}=x]
 \le(1-b_p\alpha)\|x\|_Q^p+
             C_p\alpha^2\|x\|_Q^{p-2}+C_p\alpha^p .
\]

\smallskip\noindent\emph{Iterate the moment recursion.}
For initial error $e_0=0$ put $u_{p,k}=\E\|e_k\|_Q^p$ and $u_{0,k}=1$.
The preceding inequality implies
\[
 u_{p,k}\le(1-b_p\alpha_k)u_{p,k-1}
               +C_p\alpha_k^2 u_{p-2,k-1}+C_p\alpha_k^p.
\]
If $u_{p-2,k-1}=O(\alpha_{k-1}^{(p-2)/2})$, the forcing is
$O(\alpha_k^{1+p/2})$, also for $p=2$. Moreover,
\[
 (1-b_p\alpha_k)\alpha_{k-1}^{p/2}
 \le \alpha_k^{p/2}(1-b_p\alpha_k/2)
\]
for all sufficiently large $k$ by \eqref{eq:step-relative}.
Choosing a sufficiently large multiple of $\alpha_k^{p/2}$ as a
supersolution proves $u_{p,k}\le C_p\alpha_k^{p/2}$ by induction in
even $p$, starting with $u_{0,k}=1$.

\smallskip\noindent\emph{Arbitrary initialization and the average.}
For the homogeneous recursion the bound
$\|y\|_Q\le C\alpha\|x\|_Q$ yields a multiplicative conditional
contraction $(1-b_p\alpha_k)\|x\|_Q^p$. Iteration on the vectors of
a fixed basis, with norm equivalence in fixed dimension, gives
$\|\Gamma_{i:k}\|_{L^p}\le C_p e^{-b(s_k-s_{i-1})}$.
The additional term for initial error $e_0$ is $\Gamma_{1:k}e_0$,
giving the last-iterate bound in \eqref{eq:inputs}. By independence,
$G_{i:k}=\E\Gamma_{i:k}$. The deterministic summation estimate in
Lemma~\ref{lem:kernel} now gives $\sup_{n,j}\|Q_j\|<\infty$ directly
from this product bound.

For the averaged-iterate estimate write the centered part in
\eqref{eq:mart} as $\sum_{j<n}D_{j,n}$, with
$D_{j,n}=-Q_j(\varepsilon_j+\widetilde A_j e_{j-1})/\sqrt n$.
These are martingale differences with
$\|D_{j,n}\|_{L^p}\le C_p(1+\|e_0\|)/\sqrt n$ by the last-iterate bound. For completeness
put $M=\max_{r<n}\|\sum_{j=1}^rD_{j,n}\|_{L^p}$ and
$S^2=\sum_{j<n}\|D_{j,n}\|_{L^p}^2$. Expand the Euclidean $p$th power about each partial sum and telescope.
Every linear term has zero expectation by the martingale property. H\"older's inequality on the remainders gives
\[
 M^p\le C_p M^{p-2}S^2+C_p\sum_{j<n}\|D_{j,n}\|_{L^p}^p
      \le C_p(M^{p-2}S^2+S^p).
\]
Thus $M\le C'_p S\le C''_p(1+\|e_0\|)$. Equation~\eqref{eq:mart} also gives
\[
 \|\E T_n\|\le n^{-1/2}\sum_{k=0}^{n-1}\|G_{1:k}\|\|e_0\|
             =O(n^{-1/2}).
\]
For arbitrary real $p\ge2$, apply the estimates at $q=2\lceil p/2\rceil$
and use $\|Y\|_{L^p}\le\|Y\|_{L^q}$.
For the bootstrap, the independent, bounded multipliers satisfy
$\E w_j=1$ and are independent of the observations. The pairs
$(w_j A_j,w_j\varepsilon_j)$ therefore obey the same one-step
inequalities under the joint law, with mean matrix $\bar A$.
\end{proof}

\subsection{Resampling independent coordinates}\label{app:resampling}
Let $N\ge1$, $[N]=\{1,\ldots,N\}$, and let
$X=(X_1,\ldots,X_N)$ have independent coordinates with respective laws
$\pi_i$ on measurable spaces $\mathcal X_i$. Let $X'$ be an independent copy.
For a set of indices $U\subseteq[N]$, define $X^U$ by
$(X^U)_i=X'_i$ for $i\in U$ and $(X^U)_i=X_i$ otherwise.
Thus $X^\varnothing=X$. For a function $f$ on the product space set
\begin{equation}\label{eq:resample-differences}
\begin{aligned}
 \Delta_jf(X^U)&=f(X^U)-f(X^{U\cup\{j\}}), &&j\notin U,\\
 \Delta_i\Delta_j f(X)
   &=f(X)-f(X^{\{i\}})-f(X^{\{j\}})+f(X^{\{i,j\}}), &&i\ne j .
\end{aligned}
\end{equation}
A uniform bound on differences holds for every choice of original and
replacement values in sets of full marginal probability.

Independently of $(X,X')$, draw a single index $I_*$ with
$\Pp(I_*=j)=1/N$ for each $j\in[N]$, and draw a uniform permutation
$(\tau_1,\ldots,\tau_N)$ of $[N]$, independently of $I_*$.
If $\tau_r=I_*$, define $\mathcal I=\{\tau_1,\ldots,\tau_{r-1}\}$,
which is empty if $r=1$. Conditional on $I_*=j$, its position is uniform
on $[N]$, and given the position $r$, each $(r-1)$-element subset of
$[N]\setminus\{j\}$ has the same probability. Hence
\begin{equation}\label{eq:resample-order}
 \Pp(\mathcal I=U\mid I_*=j)
  =\frac{1}{N\binom{N-1}{|U|}},\qquad U\subseteq[N]\setminus\{j\}.
\end{equation}

\begin{lemma}[Covariance under resampling]\label{lem:covid}
Suppose $f$ and $g$ are measurable real-valued functions with $\E[f(X)^2]+\E[g(X)^2]<\infty$. Then
\begin{equation}\label{eq:chatterjee-identity}
 \Cov(f(X),g(X))
 =\frac N2\E\big[\Delta_{I_*}f(X^{\mathcal I})\,\Delta_{I_*}g(X)\big].
\end{equation}
Both differences use the same replacement $X'_{I_*}$.
For vector-valued $f$ and $g$, the covariance matrix is given by the
corresponding outer product, with the second difference transposed.
\end{lemma}
The scalar identity is Lemma~2.3 of \cite{chatterjee}, proved in its
Section~4.1. See also \cite[Lemma~1.1]{kp} for the vector form.

\begin{lemma}[Efron--Stein and bounded differences]\label{lem:resamplingvariance}
For any square-integrable scalar $u(X)$,
\[
 \operatorname{Var}(u(X))\le \tfrac12\sum_{i=1}^N\E|\Delta_i u(X)|^2.
\]
If $|\Delta_i u(X)|\le c_i$ uniformly and $\sum_i c_i^2>0$, then
\[
 \Pp\{|u(X)-\E u(X)|\ge t\}\le
       2\exp\!\left(-2t^2/\sum_{i=1}^N c_i^2\right),\qquad t>0.
\]
\end{lemma}
Both estimates are proved in \cite[Theorems~3.1 and 6.2]{blm}.

\subsection{The multivariate Stein bound of Fang and R\"ollin}
All vector norms in this subsection are Euclidean. Write $\eta\sim\N(0,I_d)$.
A $d$-dimensional Stein coupling is a triple $(S,S',G)$ of
square-integrable $\mathbb R^d$-valued random vectors on one probability
space. The auxiliary vector $G$ links the change $D=S'-S$ to $S$
through the identity
\begin{equation}\label{eq:stein-coupling}
 \E[S_r g(S)]=\E[G_r\{g(S')-g(S)\}],\qquad 1\le r\le d,
\end{equation}
for scalar-valued measurable $g:\mathbb R^d\to\mathbb R$ for which
the terms are integrable. Taking $g$ bounded suffices to verify the
identity, since truncating $g$ and applying dominated convergence
then gives its integrable form.

\begin{lemma}[{Fang and R\"ollin, \cite[Theorem~2.1]{fr}}]\label{lem:stein}
Let $(S,S',G)$ be a $d$-dimensional Stein coupling with $\E S=0$,
$\Cov(S)=I_d$, $\|G\|\le a_0$ and $\|D\|\le b_0$ almost surely. Put
\begin{align*}
 B_1^2&=\operatorname{Var}(\E[\|D\|^2\mid S]),\\
 B_2^2&=\sum_{r,s=1}^d\operatorname{Var}(\E[G_rD_s\mid S]),\\
 B_3^2&=\sum_{r,s,t=1}^d\operatorname{Var}(\E[G_rD_sD_t\mid S]).
\end{align*}
With a universal constant $c$,
\begin{equation}\label{eq:fr-bound}
 \dc(S,\eta)\le c\left\{
 d^{7/4}a_0\E\|D\|^2+d^{1/4}b_0+d^{7/8}\sqrt{a_0B_1}
             +d^{3/8}B_2+d^{1/8}\sqrt{B_3}\right\}.
\end{equation}
\end{lemma}
The coupling identity is the coordinate form of
\cite[Definition~2.1]{fr}. The proof of their multivariate
bound \eqref{eq:fr-bound} is in Section~4 of \cite{fr}.

The following consequence combines the multivariate Stein bound with
the independent-coordinate covariance identity. Its hypotheses apply directly
to $f(X)$ before covariance normalization, even when the conditional
covariance depends on the observations.

\begin{lemma}
\label{lem:sensitivity}
Let $X$ have the independent coordinates of Appendix~\ref{app:resampling}.
Suppose $f(X)\in\mathbb R^d$ is centered and square integrable,
with $V=\Cov(f(X))$ and $\lambda_{\min}(V)\ge v>0$. Assume that
for some $L,H<\infty$ and nonnegative $(h_{ij})_{i,j=1}^N$,
uniformly in the replacements of \eqref{eq:resample-differences},
\begin{equation}\label{eq:uppersens}
\begin{aligned}
 \|\Delta_i f(X)\| &\le L N^{-1/2} &&(1\le i\le N),\\
 \|\Delta_i\Delta_j f(X)\| &\le L N^{-1/2}h_{ij} &&(i\ne j),\\
 h_{ii}&=1,\qquad \max_{1\le i\le N}\sum_{j=1}^Nh_{ij}\le H .
\end{aligned}
\end{equation}
Then $\dc(f(X),V^{1/2}\eta)\le C(d,v,L,H)N^{-1/2}$.
\end{lemma}
\begin{proof}
\emph{1. Construct the coupling.}
Normalize $f(X)$ by setting $f_V(x)=V^{-1/2}f(x)$. The differences in \eqref{eq:uppersens}
satisfy the same bounds for $f_V$ with $L$ replaced by $L/\sqrt v$.
With the index and set of \eqref{eq:resample-order}, take
\[
 S=f_V(X),\qquad S'=f_V(X^{\{I_*\}}),\qquad
 D=-\Delta_{I_*} f_V(X),\qquad G=-\tfrac N2\Delta_{I_*}f_V(X^{\mathcal I}).
\]
For a bounded measurable $g:\mathbb R^d\to\mathbb R$, apply
Lemma~\ref{lem:covid} to $(f_V)_r$ and $g\circ f_V$. Since
$\E S_r=0$, their covariance equals $\E[S_rg(S)]$. The two minus
signs give exactly the right-hand side of \eqref{eq:stein-coupling}.
Also $\Cov(S)=I_d$, $\|D\|\le C/\sqrt N$, $\|G\|\le C\sqrt N$ and
$\E\|D\|^2\le C/N$.

\smallskip\noindent\emph{2. Change one original coordinate.}
Fix $i\in[N]$. Take $X_i''\sim\pi_i$ independent of all variables
already defined and set $\widetilde X=(X_1,\ldots,X_{i-1},X_i'',
X_{i+1},\ldots,X_N)$. For a function $U$ of $X$ and
$\zeta=(X',I_*,\mathcal I)$ let
\[
 \delta_iU=U(X,\zeta)-U(\widetilde X,\zeta).
\]
Thus the two evaluations share the selected index, the permutation set
and all replacement values $X'$. On $I_*=j$, if $i\ne j$ then
$\delta_iD$ is a mixed difference at $i,j$. The same holds for
$\delta_iG$ when $i\notin\mathcal I$. When $i\in\mathcal I$ the two
evaluations of $G$ coincide. If $i=j$, each difference is controlled
by the first bound of \eqref{eq:uppersens}. Consequently, on $I_*=j$,
\[
 \|\delta_iD\|\le C h_{ij}/\sqrt N,\qquad
 \|\delta_iG\|\le C\sqrt N\,h_{ij}.
\]
With tildes denoting evaluation at $\widetilde X$, we have, for example,
\[
 G_rD_s-\widetilde G_r\widetilde D_s
    =(G_r-\widetilde G_r)D_s+\widetilde G_r(D_s-\widetilde D_s).
\]
Apply the same identity to three factors and use
$|\|D\|^2-\|\widetilde D\|^2|\le(\|D\|+\|\widetilde D\|)\|\delta_iD\|$.
For $1\le r,s,t\le d$, this gives
\[
 |\delta_i(\|D\|^2)|\le C h_{ij}/N,\qquad
 |\delta_i(G_rD_s)|\le C h_{ij},\qquad
 |\delta_i(G_rD_sD_t)|\le C h_{ij}/\sqrt N.
\]

\smallskip\noindent\emph{3. Average over the selected coordinate.}
For each scalar statistic $U=\|D\|^2$, $U=G_rD_s$, or $U=G_rD_sD_t$,
write $q_U(x)=\E_\zeta[U(x,\zeta)]$. Then
$q_U(X)-q_U(\widetilde X)=\E[\delta_iU\mid X,\widetilde X]$,
where integration is over the common auxiliary variables. In each
bound from Step~2, averaging over $I_*$ replaces $h_{ij}$ by
$N^{-1}\sum_{j=1}^N h_{ij}\le H/N$. The resulting bounds, uniform in
$i$ and in the fixed component indices, are
\[
\begin{array}{c|ccc}
 U&\|D\|^2&G_rD_s&G_rD_sD_t\\ \hline
 \sup_i\|q_U(X)-q_U(\widetilde X)\|_\infty
   & O(N^{-2})&O(N^{-1})&O(N^{-3/2})\\
 \operatorname{Var}(q_U(X))^{1/2}
   & O(N^{-3/2})&O(N^{-1/2})&O(N^{-1}).
\end{array}
\]
Here $\|\cdot\|_\infty$ is the essential supremum. The variance row
follows from Lemma~\ref{lem:resamplingvariance}, with the squared
coordinate bounds summed over $i$. Because $S$ is a function of $X$,
\[
 \operatorname{Var}(\E[U\mid S])
 =\operatorname{Var}(\E[q_U(X)\mid S])
 \le \operatorname{Var}(q_U(X)).
\]

\smallskip\noindent\emph{4. Apply the Stein bound.}
Summing over the component indices gives $B_1=O(N^{-3/2})$,
$B_2=O(dN^{-1/2})$ and $B_3=O(d^{3/2}N^{-1})$, with constants
depending on $v,L,H$. Inserting these bounds and
$(a_0,b_0)=(C\sqrt N,C/\sqrt N)$ in \eqref{eq:fr-bound}
gives $\dc(S,\eta)\le C(d,v,L,H)N^{-1/2}$. Applying the invertible
map $V^{1/2}$ to both laws preserves convex distance.
\end{proof}

\subsection{Gaussian comparison, translations and boundary strips}
For a convex Borel set $D\subseteq\mathbb R^d$ and $t>0$ put
\[
 D^t=\{z:\operatorname{dist}(z,D)\le t\},\qquad
 D^{-t}=\{z\in D:\operatorname{dist}(z,D^c)>t\}.
\]
Distance is Euclidean and $\operatorname{dist}(z,\varnothing)=+\infty$.
Write $Z_V=V^{1/2}\eta$ and $p_V$ for its density when $V\succ0$.
\begin{lemma}
\label{lem:gausstools}
Suppose $\nu I_d\preceq V,V'\preceq M I_d$ with $0<\nu\le M<\infty$.
There is $C=C(d,\nu,M)$ such that
\begin{align}
 \dc(Z_V,Z_{V'})&\le C\|V'-V\|,\label{eq:gauss-cov-aux}\\
 \sup_{D\in\mathrm{Conv}(\mathbb R^d)}
 \max\{\Pp(Z_V\in D^t\setminus D),
        \Pp(Z_V\in D\setminus D^{-t})\}&\le Ct\qquad(t>0).
        \label{eq:gauss-strip-aux}
\end{align}
For every Borel set $D$, $g_{D,V}(u)=\Pp(u+Z_V\in D)$ is twice
continuously differentiable and
\[
 \sup_{D,u}\{\|\nabla g_{D,V}(u)\|+\|\nabla^2 g_{D,V}(u)\|\}\le C.
\]
\end{lemma}
\begin{proof}
\emph{Covariance comparison.}
Put $H=V'-V$ and $V_s=V+sH$ for $0\le s\le1$. Then
\[
 \partial_s p_{V_s}(z)=\tfrac12p_{V_s}(z)
    \{z^\top V_s^{-1}HV_s^{-1}z-\tr(V_s^{-1}H)\}.
\]
The integral of the absolute derivative is bounded by
$2^{-1/2}\|V_s^{-1/2}HV_s^{-1/2}\|_{\rm F}$
by the variance formula for a Gaussian quadratic form, where $\|\cdot\|_{\rm F}$ is the Frobenius norm.
Integrating in $s$ gives \eqref{eq:gauss-cov-aux}.

\smallskip\noindent\emph{Boundary strips and translations.}
For $V=I_d$, the Gaussian strip estimate is
\cite[eq.~(4.10), proof of Lemma~4.2]{fr}, based on \cite{ball,bentkus}.
The map $V^{-1/2}$ sends each strip into the corresponding strip
of width $t/\sqrt\nu$ about $V^{-1/2}D$, yielding \eqref{eq:gauss-strip-aux}.
Finally $g_{D,V}(u)=\int_D p_V(z-u)\,dz$. The first and second
derivatives of $p_V$ are the density times polynomials of degree one
and two, with uniformly bounded coefficients on the stated spectral
range. Their $L^1$ norms bound the corresponding derivatives of $g_{D,V}$.
\end{proof}

\begin{lemma}
\label{lem:ordered-gauss}
Let $0\prec V\preceq V'$, and put
\[
 E_* = \{z:z^\top(V^{-1}-V'^{-1})z
                 \le\log(\det V'/\det V)\}.
\]
Then $E_*$ is convex and
\begin{align*}
 \dc(Z_V,Z_{V'})
 &=\Pp(Z_V\in E_*)-\Pp(Z_{V'}\in E_*)\\
 &=d_{\rm TV}(\law(Z_V),\law(Z_{V'}))
   =\tfrac12\int_{\mathbb R^d}|p_V(z)-p_{V'}(z)|\,dz .
\end{align*}
Here $d_{\rm TV}(\mu,\nu)=\sup_D|\mu(D)-\nu(D)|$, with $D$ ranging
over Borel sets. See \cite[Proposition~6.15]{wu} for strict covariance
order.
\end{lemma}
\begin{proof}
Taking the logarithm of the density ratio gives
$\{p_V\ge p_{V'}\}=E_*$. The matrix $V^{-1}-V'^{-1}$ is positive
semidefinite, so this sublevel set is convex. Since the two densities
both integrate to one, integrating $p_V-p_{V'}$ over $E_*$ gives
$\tfrac12\int|p_V-p_{V'}|$. This attains the supremum over all Borel
sets, and $E_*$ is among the admissible convex sets.
\end{proof}

\section{Technical estimates for the LSA recursion}\label{app:finite}
Let $e_k=\theta_k-\theta^\star$, $\mathcal F_0=\{\varnothing,\Omega\}$,
$\mathcal F_k=\sigma(Z_1,\ldots,Z_k)$, and $s_k=\sum_{j=1}^k\alpha_j$
with $s_0=0$. We use the estimates in \eqref{eq:inputs}.
The constants below are uniform in the time indices and in $n$.

\subsection{Two sums of deterministic transition weights}
\begin{lemma}
\label{lem:kernel}
For each $b>0$ and the step sizes of A3, both
\[
 \sup_{k\ge1}\sum_{j=1}^k\alpha_j e^{-b(s_k-s_j)}
 \quad\text{and}\quad
 \sup_{i\ge1}\alpha_i\sum_{k=i}^\infty e^{-b(s_k-s_i)}
\]
are finite. Consequently $\sup_{n\ge2}\max_{1\le j<n}\|Q_j\|<\infty$.
\end{lemma}
\begin{proof}
Since $0<\alpha_j\le\alpha_1$, there is a $c_b>0$ such that
$1-e^{-b\alpha_j}\ge c_b\alpha_j$. Hence
\[
 c_b\sum_{j=1}^k\alpha_j e^{-b(s_k-s_j)}
 \le\sum_{j=1}^k\left[e^{-b(s_k-s_j)}-e^{-b(s_k-s_{j-1})}\right]
 =1-e^{-b s_k}\le1 .
\]
For $i\le k\le2i$, $s_k-s_i\ge(k-i)\alpha_{2i}$ and
$\alpha_{2i}\ge2^{-\gamma}\alpha_i$. Thus
\[
 \sum_{k=i}^{2i} e^{-b(s_k-s_i)}
 \le (1-e^{-b\alpha_{2i}})^{-1}\le C_b/\alpha_i .
\]
For $k>2i$, monotonicity gives $s_k-s_i\ge(k-i)\alpha_k\ge k\alpha_k/2$.
The series $\sum_{k\ge1}\exp(-bk\alpha_k/2)$ converges because
$k\alpha_k\asymp k^{1-\gamma}$. This bounds the remaining tail
uniformly in $i$. Combining the two estimates and applying
\eqref{eq:inputs} to
$\|Q_j\|\le\alpha_j\sum_{k=j}^{n-1}\|G_{j+1:k}\|$ proves the result.
\end{proof}

\subsection{Covariance of the normalized average}
Let $Y_j=\varepsilon_j+\widetilde A_j e_{j-1}$, $\mu_k=\E e_k$ and
$M_k=\E[e_ke_k^\top]$.
For a deterministic vector $u$ and matrix $M$, define the linear maps
\[
 \mathcal B(u)=\E[\varepsilon_1 u^\top\widetilde A_1^\top],
 \qquad
 \mathcal V(M)=\E[\widetilde A_1M\widetilde A_1^\top].
\]
By A2, $\|\mathcal B(u)\|\le C\|u\|$ and $\|\mathcal V(M)\|\le C\|M\|$.
Independence of the new observation from the past and $\mu_k=G_{1:k}e_0$ give
\[
 \E[Y_jY_j^\top]=
 \Sigma_\varepsilon+\mathcal V(M_{j-1})
 +\mathcal B(\mu_{j-1})+\mathcal B(\mu_{j-1})^\top.
\]
In particular $\E[\varepsilon_j e_{j-1}^\top\widetilde A_j^\top]
=\mathcal B(\mu_{j-1})$, because only the past factor $e_{j-1}$
is independent of the present pair. For $i<j$, conditioning on
$\mathcal F_{j-1}$ gives $\E[Y_iY_j^\top]=0$. It follows that
\begin{equation}\label{eq:covexact}
 \Cov(T_n)-\Sigma_n
 =\frac1n\sum_{j=1}^{n-1}Q_j
 \bigl[\mathcal V(M_{j-1})+\mathcal B(\mu_{j-1})
              +\mathcal B(\mu_{j-1})^\top\bigr] Q_j^\top .
\end{equation}
Here $\sum_{j<n}\|\mu_{j-1}\|\le C$ by product decay. Also
$\sum_{j<n}\|M_{j-1}\|\le C(1+\sum_{k=1}^{n-1}\alpha_k)\le Cn^{1-\gamma}$
by \eqref{eq:inputs}. The uniform bound on $Q_j$ now yields
$\|\Cov(T_n)-\Sigma_n\|\le C(n^{-1}+n^{-\gamma})$.
In particular, at $\gamma=1/2$ the quadratic term is controlled by
$n^{-1}\sum_{k=1}^{n-1}\alpha_k\sim 2c_0 n^{-1/2}$.

\subsection{The restarted recursion and its exact remainder}\label{sec:restartexp}
Fix $0\le m\le n-2$ and $x\in\mathbb R^d$.
The observations after time $m$ form the vector
$X=(Z_{m+1},\ldots,Z_{n-1})$ with $N=n-1-m$ independent coordinates.
Write $\Pp_x$ and $\E_x$ for probability and expectation over
$X\sim\pi^{\otimes N}$ with $x$ held fixed. Define the restarted errors
$e_k^{m,x}$, one vector for each $m\le k<n$, by
\begin{equation}\label{eq:restartpath}
 e_m^{m,x}=x,\qquad
 e_k^{m,x}=(I_d-\alpha_k A_k)e_{k-1}^{m,x}-\alpha_k\varepsilon_k
 \quad(m<k<n).
\end{equation}
Taking $x=e_m$ and using the observed inputs after $m$ gives $e_k^{m,e_m}=e_k$
pathwise. Thus conditional statements about the original recursion given
$\mathcal F_m$ can be obtained from deterministic-start statements at $x=e_m$.

\paragraph{Derivation of the remainder.}
For fixed $m,x$, abbreviate $J_k^{(\ell)}=J_{m,k}^{(\ell)}(x)$,
$H_k^{(K)}=H_{m,k}^{(K)}(x)$ and $B_k=I_d-\alpha_k\bar A$.
Let $S_k^{[K]}=\sum_{\ell=0}^KJ_k^{(\ell)}$. The initial conditions
and updates are those of \eqref{eq:J}--\eqref{eq:H}.
To determine the transition for $H^{(K)}$, first add the equations in
\eqref{eq:J}. Collecting the terms containing $\widetilde A_k$ gives
\[
 S_k^{[K]}
  =(I_d-\alpha_k A_k)S_{k-1}^{[K]}-\alpha_k\varepsilon_k
      +\alpha_k\widetilde A_k J_{k-1}^{(K)} .
\]
Subtracting this equality from \eqref{eq:restartpath} and iterating yields
\begin{equation}\label{eq:restartH}
\begin{aligned}
 H_m^{(K)} &=0,\\
 H_k^{(K)} &=(I_d-\alpha_k A_k)H_{k-1}^{(K)}
                    -\alpha_k\widetilde A_k J_{k-1}^{(K)}
     =-\sum_{j=m+1}^k\alpha_j\Gamma_{j+1:k}\widetilde A_jJ_{j-1}^{(K)} .
\end{aligned}
\end{equation}
In particular it is the random transition $I_d-\alpha_k A_k$ that propagates
the remainder. Pathwise, at every time and at every truncation order,
\[
 e_k^{m,x}=\sum_{\ell=0}^K J_{m,k}^{(\ell)}(x)+H_{m,k}^{(K)}(x),
 \qquad
 H_{m,k}^{(K)}(x)=J_{m,k}^{(K+1)}(x)+H_{m,k}^{(K+1)}(x).
\]
At $m=0,x=0$ these are precisely the fluctuation terms in
\cite[eqs.~(3.11)--(3.13) and (5.1)--(5.2)]{butyrin}. With $x\ne0$, the
zeroth term also includes $G_{m+1:k}x$ and subsequent terms propagate its
interaction with the centered matrices.

On summing over $m\le k<n$ and dividing by $\sqrt n$, we obtain
\eqref{eq:sketch-FH}. Its $k=m$ term contributes $x/\sqrt n$.
All these sums are affine functions of $x$ for fixed observations.

\paragraph{Connection with the operator notation.}
For a sequence $v=(v_j)_{j=m}^{n-1}$ with $v_{j-1}$ measurable with
respect to $\sigma(Z_{m+1},\ldots,Z_{j-1})$, set
\[
 (\mathcal T_m v)_m=0,\qquad
 (\mathcal T_m v)_k=-\sum_{j=m+1}^k\alpha_jG_{j+1:k}
                                      \widetilde A_jv_{j-1}.
\]
Unrolling each mean-matrix recursion in \eqref{eq:J} gives
\[
 J_{m,k}^{(0)}(x)=G_{m+1:k}x-\sum_{j=m+1}^k\alpha_jG_{j+1:k}\varepsilon_j,
 \qquad J_{m,k}^{(\ell+1)}(x)=(\mathcal T_m J_{m,\bullet}^{(\ell)}(x))_k .
\]
Here $\bullet$ denotes the full indexed sequence $m,\ldots,n-1$ and the
subscript $k$ selects its $k$th time entry. Unrolling \eqref{eq:restartpath}
around $\bar A$ similarly gives
$e_k^{m,x}=J_{m,k}^{(0)}(x)+(\mathcal T_m e_\bullet^{m,x})_k$.
Successive substitution gives the coordinate identities
\begin{equation}\label{eq:operatorbridge}
 J_{m,k}^{(\ell)}(x)=(\mathcal T_m^\ell J_{m,\bullet}^{(0)}(x))_k,\qquad
 H_{m,k}^{(K)}(x)=(\mathcal T_m^{K+1}e_\bullet^{m,x})_k .
\end{equation}
The zeroth power is the identity. We use the random-product representation
\eqref{eq:restartH} for the pathwise expansion and \eqref{eq:operatorbridge}
for the moment estimate below.

\begin{lemma}
\label{lem:proxy}
Assume A1--A3. Let $1\le m\le n-2$, $\|x\|\le1$ and fix an integer $K\ge0$. Uniformly in $m,n,x$,
\begin{equation}\label{eq:uppertrunc}
 \E_x F_x^{[K]}=\E_x F_x,\qquad
 \|F_x-F_x^{[K]}\|_{L^2(\Pp_x)}\le C_K\sqrt n\,\alpha_m^{(K+1)/2}.
\end{equation}
For replacements as in \eqref{eq:resample-differences} at times $i,j\in\{m+1,\ldots,n-1\}$, we have
\begin{align*}
 \|\Delta_i F_x^{[K]}\|&\le C_K n^{-1/2},&
 \|\Delta_i\Delta_j F_x^{[K]}\|&\le C_K n^{-1/2}h_{ij}\quad(i\ne j),\\
 h_{ii}&=1,& h_{ij}=h_{ji}&=\alpha_i e^{-b(s_{j-1}-s_i)}
                                  \quad(i<j)
\end{align*}
for some $b>0$, with $\sup_i\sum_j h_{ij}\le C_K$.
\end{lemma}
\begin{proof}
\emph{1. Moment of the exact remainder.}
For each fixed $k$ and adapted $v$, the vectors
$G_{j+1:k}\widetilde A_j v_{j-1}$, $m<j\le k$, are martingale differences
in $j$. Indeed $G_{j+1:k}$ is deterministic and
$\E[\widetilde A_j\mid Z_{m+1},\ldots,Z_{j-1}]=0$. Put
$\|v\|_{2,*}=\max_{m\le k<n}(\E_x\|v_k\|^2)^{1/2}$. Orthogonality at
distinct final input times gives
\begin{align*}
 \E_x\|(\mathcal T_m v)_k\|^2
 &=\sum_{j=m+1}^k\alpha_j^2
       \E_x\|G_{j+1:k}\widetilde A_jv_{j-1}\|^2\\
 &\le C\|v\|_{2,*}^2\alpha_m
     \sum_{j=m+1}^k\alpha_j e^{-b(s_k-s_j)}
 \le C'\alpha_m\|v\|_{2,*}^2 .
\end{align*}
The last two inequalities use decreasing steps, \eqref{eq:inputs}, and
Lemma~\ref{lem:kernel}. The unrolled recursion with random products also gives
\[
 (\E_x\|e_k^{m,x}\|^2)^{1/2}
 \le\|\Gamma_{m+1:k}\|_{L^2}\|x\|
   +\|\varepsilon\|_\infty\sum_{j=m+1}^k
                         \alpha_j\|\Gamma_{j+1:k}\|_{L^2}
 \le C .
\]
By \eqref{eq:operatorbridge},
$\|H_{m,\bullet}^{(K)}(x)\|_{2,*}\le C_K\alpha_m^{(K+1)/2}$.
Minkowski's inequality in \eqref{eq:sketch-FH} now multiplies this
bound by $(n-m)/\sqrt n\le\sqrt n$. Moreover each coordinate of
$\mathcal T_m v$ has mean zero for adapted $v$, so each $H_{m,k}^{(K)}(x)$
has mean zero. This proves both assertions of \eqref{eq:uppertrunc}.

\emph{2. Propagation of a perturbation through the retained terms.}
Collect the retained terms in a vector and equip the product space with
the maximum norm:
\[
 Y_k=(J_{m,k}^{(0)}(x),\ldots,J_{m,k}^{(K)}(x)),\qquad
 \|y\|_*=\max_{0\le r\le K}\|y_r\|.
\]
Equation~\eqref{eq:J} is the affine recursion $Y_k=M_kY_{k-1}+q_k$,
where $q_k=(-\alpha_k\varepsilon_k,0,\ldots,0)$ and
\[
 (M_k y)_0=B_ky_0,\qquad
 (M_k y)_r=B_ky_r-\alpha_k\widetilde A_ky_{r-1},
                         \quad 1\le r\le K .
\]
Write $M_{i+1:k}=M_k\cdots M_{i+1}$ with empty product the identity.
In expanding a block of this product, each subdiagonal factor raises the
level by one. Hence each term uses at most $K$ such factors. A choice of
times $i<l_1<\cdots<l_r\le k$ for these factors has norm at most
\[
 C_K e^{-c(s_k-s_i)}\prod_{a=1}^r\alpha_{l_a}.
\]
Indeed, each intervening diagonal interval is a product of mean transitions.
The factors omitted from their total exponential decay cost at most
$\exp(c\sum_{a=1}^r\alpha_{l_a})\le\exp(cK\alpha_1)$. Summing in the
positions, using
$\sum_{l_1<\cdots<l_r}\prod_a\alpha_{l_a}\le(s_k-s_i)^r/r!$ for these times,
gives the induced operator-norm bound
\[
 \|M_{i+1:k}\|_*
 \le C_K e^{-c(s_k-s_i)}\sum_{r=0}^K\frac{(s_k-s_i)^r}{r!}
 \le C_K' e^{-b(s_k-s_i)}
\]
with $0<b<c$. Applying this bound to
$Y_k=M_{m+1:k}Y_m+\sum_{j=m+1}^k M_{j+1:k}q_j$ and using
Lemma~\ref{lem:kernel} yields $\max_{m\le k<n}\|Y_k\|_*\le C_K$ pathwise.

\emph{3. One and two observation replacements.}
Use primes for the coefficients evaluated at a replaced input.
Changing input $i$ leaves $Y_{i-1}$ fixed. The difference at $i$ and its
later values are therefore
\[
 d_i=(M_i-M_i')Y_{i-1}+q_i-q_i',\qquad
 \|d_i\|_*\le C_K\alpha_i,\qquad
 \Delta_iY_k=M_{i+1:k}d_i\quad(k\ge i).
\]
For $i<j$, a second replacement at $j$ acts on this propagated difference.
The two additive changes at $j$ cancel, giving the exact mixed difference
\[
 \Delta_i\Delta_jY_k
   =M_{j+1:k}(M_j-M_j')M_{i+1:j-1}d_i,\qquad k\ge j .
\]
Consequently the bounds for the first and mixed differences are, respectively,
\[
 C_K\alpha_i e^{-b(s_k-s_i)},\qquad
 C_K\alpha_i\alpha_j e^{-b(s_{j-1}-s_i)}e^{-b(s_k-s_j)} .
\]
In \eqref{eq:sketch-FH} the map $y\mapsto\sum_{r=0}^K y_r$ has norm
at most $K+1$. Summing over $k\ge i$ for the first difference and $k\ge j$
for the mixed difference costs at most $C/\alpha_i$ and $C/\alpha_j$,
respectively, by Lemma~\ref{lem:kernel}. Division by $\sqrt n$ proves
the claimed estimates. For the row sums, at a fixed $i$ the two parts satisfy
\[
 \sum_{j<i}\alpha_j e^{-b(s_{i-1}-s_j)}\le C,\qquad
 \alpha_i\sum_{j>i}e^{-b(s_{j-1}-s_i)}\le C .
\]
These are the two assertions of Lemma~\ref{lem:kernel}.
\end{proof}

\subsection{Covariance of the restarted sum}
\begin{lemma}
\label{lem:suffixcov}
Under the hypotheses of Theorem~\ref{thm:finite}, let $1\le m<n-1$
and let $x\in\mathbb R^d$ be deterministic. For the restarted sum
$F_x=n^{-1/2}\sum_{k=m}^{n-1} e_k^{m,x}$, let
$V_{m,x}=\operatorname{Cov}_x(F_x)$, where the covariance is taken over
$Z_{m+1},\ldots,Z_{n-1}$. Then
\[
 \|V_{m,x}-\Sigma_n\|\le
 C\left(\frac mn+n^{-\gamma}+\frac{1+\|x\|^2}{n\alpha_m}\right).
\]
\end{lemma}
\begin{proof}
\emph{Centering and deleting the initial weights.}
For $j>m$, set $Y_j^{m,x}=\varepsilon_j+\widetilde A_j e_{j-1}^{m,x}$.
The same change of summation order as in \eqref{eq:mart} yields
\[
 \E_x F_x=\frac1{\sqrt n}\sum_{k=m}^{n-1}G_{m+1:k}x,\qquad
 F_x-\E_x F_x=-\frac1{\sqrt n}\sum_{j=m+1}^{n-1}Q_jY_j^{m,x}.
\]
These increments are centered conditional on the restarted past, so
\[
 V_{m,x}-\Sigma_n=
 -\frac1n\sum_{j=1}^mQ_j\Sigma_\varepsilon Q_j^\top
 +\frac1n\sum_{j=m+1}^{n-1}Q_j
       \{\E_x[Y_j^{m,x}(Y_j^{m,x})^\top]-\Sigma_\varepsilon\}Q_j^\top .
\]
The first term has norm bounded by $Cm/n$.

\emph{First and second moments on the remaining path.}
Use the maps $\mathcal B,\mathcal V$ of \eqref{eq:covexact}.
Since $\E_x e_k^{m,x}=G_{m+1:k}x$, the two cross terms in braces,
after summation and division by $n$, have combined norm at most
\[
 \frac Cn\sum_{k=m}^{n-2}\|G_{m+1:k}x\|
 \le \frac{C\|x\|}{n\alpha_m}
 \qquad\text{by Lemma~\ref{lem:kernel}.}
\]
For the quadratic term, couple the restarted path to an ordinary
trajectory with the same observations after $m$. Direct subtraction
of the two recursions gives
$e_k^{m,x}-e_k=\Gamma_{m+1:k}(x-e_m)$. The suffix product is independent
of $e_m$. The joint expectation on this coupling therefore satisfies
\[
 \begin{aligned}
 \E_x\|e_k^{m,x}\|^2
 &\le 2\E\|e_k\|^2+
      2\E\|\Gamma_{m+1:k}\|^2\,\E\|x-e_m\|^2\\
 &\le C\{\alpha_k+(1+\|x\|^2)e^{-b(s_k-s_m)}\},\qquad k\ge m .
 \end{aligned}
\]
Here \eqref{eq:inputs} absorbs the exponentially small initial error into
$C\alpha_k$. The norm of the average of the quadratic terms is at most
\[
 \frac Cn\sum_{k=m}^{n-2}\E_x\|e_k^{m,x}\|^2
 \le C n^{-\gamma}+\frac{C(1+\|x\|^2)}{n\alpha_m}.
\]
Adding all three contributions and using $\|x\|\le1+\|x\|^2$ proves
the lemma.
\end{proof}

\section{Covariance expansions and smooth approximation}
\label{app:endpoint}

\subsection{The deterministic covariance correction}
\begin{lemma}
\label{lem:matrix}
Suppose $-\bar A$ is Hurwitz, $\Sigma_\varepsilon\succ0$, and
$\alpha_k=c_0(k+k_0)^{-\gamma}$ with $c_0>0$, fixed $k_0\ge0$ and
$1/2\le\gamma<1$. For the weights in \eqref{eq:targets}, the limit in
operator norm is
\[
 n\alpha_n(\Sigma_n-\Sigma_\infty)
 \longrightarrow S:=\int_0^\infty
 e^{-t\bar A}\Sigma_\infty e^{-t\bar A^\top}\,dt .
\]
In particular, $S=c_0L$ with $L$ given by \eqref{eq:lyap}, and
\eqref{eq:exp} holds.
\end{lemma}
\begin{proof}
Wu et al.\ \cite[Theorem~6.11]{wu} give the coefficient for their
TD averaging convention. We verify it for the finite sum defining $Q_j$
in this paper.

\smallskip\noindent\emph{Cancellation of the linear terms.}
Put $D_j=Q_j-\bar A^{-1}$ and
$s_k=\sum_{i=1}^k\alpha_i$. The factors in $G$ commute with $\bar A$.
Telescoping in $j$, followed by interchanging the two sums, gives
\begin{align*}
 \alpha_jG_{j+1:k}\bar A &=G_{j+1:k}-G_{j:k},\\
 \sum_{j=1}^{n-1}Q_j
 &=\sum_{k=1}^{n-1}\sum_{j=1}^k\alpha_jG_{j+1:k}
   =\sum_{k=1}^{n-1}(I-G_{1:k})\bar A^{-1} .
\end{align*}
The deterministic product bound $\|G_{j+1:k}\|\le
C\exp[-b(s_k-s_j)]$ follows from Lyapunov contraction for sufficiently
small steps, with a finite initial product absorbed in $C$. In
particular $\sum_{k\ge1}\|G_{1:k}\|<\infty$. Thus
$\sum_{j<n}D_j=-\sum_{k<n,\,k\ge1}G_{1:k}\bar A^{-1}=O(1)$. Expanding
$Q_j=\bar A^{-1}+D_j$ in the covariance gives
\begin{align}
 n(\Sigma_n-\Sigma_\infty)
 &=\sum_{j=1}^{n-1} D_j\Sigma_\varepsilon D_j^\top+O(1).
 \label{eq:lower-squares}
\end{align}
Here the remainder is exactly the sum of
$(\sum_{j<n}D_j)\Sigma_\varepsilon \bar A^{-\top}$, its transpose, and
$-\Sigma_\infty$.

\smallskip\noindent\emph{Contribution of the final iterates.}
Telescoping instead in $k$ gives
$\sum_{k=j}^{n-1}\alpha_{k+1}G_{j+1:k}=(I-G_{j+1:n})\bar A^{-1}$.
Subtracting this identity from the definition of $Q_j$ yields
\[
 D_j=-G_{j+1:n}\bar A^{-1}+R_j,\qquad
 R_j=\sum_{k=j}^{n-1}(\alpha_j-\alpha_{k+1})G_{j+1:k}.
\]
The mean value theorem gives
$0\le\alpha_j-\alpha_{k+1}\le C\alpha_j(k-j+1)/j$.
For $j\le k\le2j$, $s_k-s_j\ge(k-j)\alpha_{2j}$ and
$\alpha_{2j}\asymp\alpha_j$. For $k>2j$, $s_k-s_j\ge c k^{1-\gamma}$
once $k$ is large. Therefore
\[
 \sum_{k\ge j}(k-j+1)e^{-b(s_k-s_j)}\le C\alpha_j^{-2},
 \qquad
 \|R_j\|\le C/(j\alpha_j).
\]
In particular $\sum_{j<n}\|R_j\|^2$ is $O(\log n)$ at $\gamma=1/2$
and $O(n^{2\gamma-1})$ for $1/2<\gamma<1$. Both bounds are
$o(\alpha_n^{-1})$. Also $s_n-s_j\ge(n-j)\alpha_n$, so
$\sum_{j<n}\|G_{j+1:n}\|^2=O(\alpha_n^{-1})$. By Cauchy--Schwarz,
\[
 \alpha_n\sum_{j<n}\|G_{j+1:n}\|\|R_j\|
 \le
 \left(\alpha_n\sum_{j<n}\|G_{j+1:n}\|^2\right)^{1/2}
 \left(\alpha_n\sum_{j<n}\|R_j\|^2\right)^{1/2}=o(1).
\]
Using \eqref{eq:lower-squares}, we obtain
\[
 n\alpha_n(\Sigma_n-\Sigma_\infty)
  =\alpha_n\sum_{j=1}^{n-1}
       G_{j+1:n}\Sigma_\infty G_{j+1:n}^\top+o(1).
\]

\smallskip\noindent\emph{Passage to the Lyapunov integral.}
Set $r=n-j$ and first restrict to $r\alpha_n\le M$ for fixed $M$.
Here $r/n=o(1)$ and $\alpha_i/\alpha_n\to1$ uniformly over
$n-r<i\le n$. Consequently
\[
 \sum_{i=n-r+1}^n\alpha_i=r\alpha_n+o(1),\qquad
 \sum_{i=n-r+1}^n\alpha_i^2=O_M(\alpha_n).
\]
Since $\log(I-\alpha_i\bar A)=-\alpha_i\bar A+O(\alpha_i^2)$ and
all factors commute, $G_{n-r+1:n}=e^{-r\alpha_n\bar A}+o(1)$ uniformly
on this range. Its contribution is a Riemann sum over $[0,M]$ with mesh
$\alpha_n$. The norm of the remaining part is at most
$C\alpha_n\sum_{r>M/\alpha_n} e^{-2br\alpha_n}\le C'e^{-2bM}$.
Letting $n\to\infty$ and then $M\to\infty$ gives $S$.

The integrand is positive definite at $t=0$ and continuous, so
$S\succ0$. For $K(t)=e^{-t\bar A}\Sigma_\infty e^{-t\bar A^\top}$,
$K'(t)=-\bar A K(t)-K(t)\bar A^\top$. Integration gives
$\bar A S+S\bar A^\top=\Sigma_\infty$. Equation~\eqref{eq:lyap}
and $n^\gamma\alpha_n\to c_0$ now give $S=c_0L$ and \eqref{eq:exp}.
\end{proof}

\subsection{Smooth approximation at the endpoint}
\begin{lemma}\label{lem:endpoint-smooth}
Under A1--A4 with $\gamma=1/2$, fix a unit vector $v$.
Let $m_n=\E[v^\top T_n]$, $\tau_n^2=\operatorname{Var}(v^\top T_n)$
and $Z\sim\N(0,1)$. There is a constant $C$ independent of $f$ and $n$
such that every bounded $C^3$ function $f$ with bounded derivatives satisfies
\[
 \big|\E f(v^\top T_n)-\E f(m_n+\tau_n Z)\big|
 \le C n^{-1/2}\|f'''\|_\infty .
\]
\end{lemma}
\begin{proof}
\emph{Replace one martingale increment at a time.}
Put $N=n-1$, $q_j=Q_j^\top v$,
$X_j=-q_j^\top Y_j/\sqrt n$ with $Y_j$ from \eqref{eq:mart},
and $S_k=\sum_{i=1}^kX_i$, $S_0=0$. Thus $v^\top T_n=m_n+S_N$.
The martingale differences are orthogonal, so
$\tau_n^2=\sum_{j=1}^N\E X_j^2$. By Lemmas~\ref{lem:moments}
and~\ref{lem:kernel}, $\sup_j\|q_j\|\le C$,
$\sup_j\|X_j\|_{L^3}\le C/\sqrt n$, and
$\sup_k\|S_k\|_{L^2}\le C$.
Take independent $\zeta_j\sim\N(0,\E X_j^2)$, independent of the observations,
and write
\[
 f_j(x)=\E f\!\left(m_n+x+\sum_{i=j+1}^{N}\zeta_i\right),
 \qquad D_j=f_j'',\qquad K_f=\|f'''\|_\infty.
\]
Telescoping over the replacements gives
\[
 \E f(v^\top T_n)-\E f(m_n+\tau_n Z)
 =\sum_{j=1}^{N}\E\big[f_j(S_{j-1}+X_j)-f_j(S_{j-1}+\zeta_j)\big].
\]
The linear terms in the Taylor expansions have zero expectation.
The third-order remainders sum to at most
$C K_f\sum_j(\E|X_j|^3+\E|\zeta_j|^3)\le C K_f/\sqrt n$.
It remains to bound
\begin{equation}\label{eq:endpoint-conditional-variance}
 \frac12\sum_{j=1}^{N}\E\!\left[D_j(S_{j-1})
    \big\{\E[X_j^2\mid\mathcal F_{j-1}]-\E X_j^2\big\}\right].
\end{equation}

\smallskip\noindent\emph{The quadratic term in the past iterate.}
For a deterministic vector $x$, define
\[
 l_j(x)=2\E[(q_j^\top\varepsilon_1)(q_j^\top\widetilde A_1x)],
 \qquad H_j(x)=\E[(q_j^\top\widetilde A_1x)^2].
\]
Then $n\E[X_j^2\mid\mathcal F_{j-1}]$ is
$q_j^\top\Sigma_\varepsilon q_j+l_j(e_{j-1})+H_j(e_{j-1})$,
where $|l_j(x)|\le C\|x\|$ and $|H_j(x)|\le C\|x\|^2$.
The constant term disappears under centering. Since $D_j$ is
$K_f$-Lipschitz, subtracting $D_j(0)$ and using Cauchy--Schwarz gives
\[
 |\operatorname{Cov}(D_j(S_{j-1}),H_j(e_{j-1}))|
 \le C K_f\|S_{j-1}\|_{L^2}\|e_{j-1}\|_{L^4}^2 .
\]
Lemma~\ref{lem:moments} gives $\sum_{j=1}^{N}\|e_{j-1}\|_{L^4}^2
\le C(1+\sum_{j=1}^{N}\alpha_j)\le C\sqrt n$.

\smallskip\noindent\emph{The linear term in the past iterate.}
For $i<j$, let $S_{j-1}^{(i)}$ be computed after replacing observation
$i$ by an independent copy. Denote by $\Delta_i$ the original quantity
minus its value after this replacement. At time $i$,
$\Delta_i e_i=-\alpha_i\Delta_iY_i$, with $\|\Delta_iY_i\|_{L^2}\le C$.
For $i<r<j$ the recursion gives
\[
 \Delta_iY_r=\widetilde A_r\Gamma_{i+1:r-1}(-\alpha_i\Delta_iY_i).
\]
The future product is independent of $\Delta_iY_i$. Lemmas~\ref{lem:moments}
and~\ref{lem:kernel} therefore imply
\[
 \|S_{j-1}-S_{j-1}^{(i)}\|_{L^2}
 \le \frac C{\sqrt n}
 \left(1+\alpha_i\sum_{r=i+1}^{j-1}e^{-b(s_{r-1}-s_i)}\right)
 \le \frac{C'}{\sqrt n}.
\]
Conditionally on all observations except $Z_i$ and on its independent
replacement, $Y_i$ has mean zero.
Consequently,
\[
 \big\|\E[Y_i D_j(S_{j-1})]\big\|
 =\big\|\E[Y_i\{D_j(S_{j-1})-D_j(S_{j-1}^{(i)})\}]\big\|
 \le C K_f/\sqrt n.
\]
Iterating the recursion with the mean matrix gives
$e_{j-1}-\E e_{j-1}=-\sum_{i<j}\alpha_i G_{i+1:j-1}Y_i$.
As $l_j$ is linear, the preceding bound and
$\sum_{i<j}\alpha_i\|G_{i+1:j-1}\|\le C$ show that
\[
 |\operatorname{Cov}(D_j(S_{j-1}),l_j(e_{j-1}))|
 \le C K_f/\sqrt n .
\]
After division by $n$, the quadratic and linear contributions to
\eqref{eq:endpoint-conditional-variance} are each bounded by
$C K_f/\sqrt n$. Adding the Taylor remainders proves the claim.
\end{proof}

\subsection{Variance of the multiplicative correction in Proposition~\ref{prop:sharp-example}}
Here $\operatorname{Var}(A_j-1)=1/4$ and $A_j$ is independent of
$\xi_j$. Writing $v_k=\E\theta_k^2$, the martingale representation
\eqref{eq:mart} and the one-step recursion give
\begin{align*}
 \operatorname{Var}(T_n)-\Sigma_n
    &=\frac1{4n}\sum_{j=1}^{n-1}Q_j^2 v_{j-1},\\
 v_k&=(1-2\alpha_k+\tfrac54\alpha_k^2)v_{k-1}+\alpha_k^2,\qquad v_0=0.
\end{align*}
For $k\ge2$, division by $\alpha_k$ gives
\[
 \frac{v_k}{\alpha_k}
  =\big(1-2\alpha_k+o(\alpha_k)\big)
         \frac{v_{k-1}}{\alpha_{k-1}}+\alpha_k,
\]
since $\alpha_{k-1}/\alpha_k=1+O(k^{-1})=1+o(\alpha_k)$.
By \cite[Lemma~B.2]{bns}, $v_k/\alpha_k\to1/2$.
On $\delta n\le j\le(1-\delta)n$ for fixed $0<\delta<1/2$,
the estimates for $D_j$ in Lemma~\ref{lem:matrix} imply $Q_j\to1$
uniformly. The weights are bounded on the full range and
$v_{j-1}\le C j^{-\gamma}$ for $j\ge2$. Riemann summation on the
interior, followed by $\delta\downarrow0$, gives
\[
 n^\gamma\{\operatorname{Var}(T_n)-\Sigma_n\}
 \longrightarrow \frac{c_0}{8}\int_0^1 t^{-\gamma}\,dt
                 =\frac{c_0}{8(1-\gamma)} .
\]

\section{Proofs of the bootstrap results}\label{app:app-proofs}

\subsection{Proof of Theorem~\ref{thm:boot}}

Put $\mathcal Z_n=\sigma(Z_1,\ldots,Z_{n-1})$ and write
$e_k=\theta_k-\theta^\star$ and $e_k^b=\theta_k^b-\theta^\star$.
For the bootstrap proof use
\[
 m=\lfloor n^{1/3}\rfloor,\qquad I_+=\{m+1,\ldots,n-1\},\qquad
 N=|I_+|=n-1-m,\qquad K=23 .
\]
Under the joint law of data and multipliers the weighted recursion has mean
matrix $\bar A$ and innovation covariance
$\E[w_j^2\varepsilon_j\varepsilon_j^\top]=2\Sigma_\varepsilon$.
In particular its averaged linear covariance in \eqref{eq:targets}
is $2\Sigma_n$. Boundedness of the weights permits application of
Lemma~\ref{lem:moments} at the fixed orders used below throughout A3.

\paragraph{1. Approximation for a deterministic restart.}
For deterministic $x\in\mathbb R^d$, use $e_k^{m,x,b}$ for the
restarted errors driven by the weighted observations, and $e_k^{m,x}$
for the restart driven by the original observations. Define the three sums
\[
 F_x^b=\frac1{\sqrt n}\sum_{k=m}^{n-1}e_k^{m,x,b},\qquad
 F_x^0=\frac1{\sqrt n}\sum_{k=m}^{n-1}e_k^{m,x},\qquad
 P_x=\frac1{\sqrt n}\sum_{k=m}^{n-1}\sum_{\ell=0}^{K}
                  J_{m,k}^{(\ell),b}(x).
\]
Here the superscript $b$ in $J$ and $H$ means that \eqref{eq:sketch-FH}
uses the inputs $(w_jA_j,w_j\varepsilon_j)$ and centered matrix
$w_j A_j-\bar A$. Thus $P_x$ is the weighted $F_x^{[K]}$ of that formula and
\[
 R_x:=F_x^b-P_x=\frac1{\sqrt n}\sum_{k=m}^{n-1}H_{m,k}^{(K),b}(x).
\]

\begin{lemma}
\label{lem:bootcov}
Let $\mathcal Z_+=\sigma(Z_j:j\in I_+)$,
$V_x=\Cov(P_x\mid\mathcal Z_+)$ and
$H_n=\sup_{\|x\|\le1}\|V_x-\Sigma_n\|$. For all sufficiently large $n$,
\begin{align}
 \sup_{\|x\|\le1}\|R_x\|_{L^2}&\le Cn^{-3/2},\label{eq:boot-trunc}\\
 \E H_n&\le Cn^{-1/2},\qquad
 \Pp\{H_n>C\sqrt{\log n/n}\}\le Cn^{-2}.
 \label{eq:bootcov-rev}
\end{align}
The $L^2$ norm in \eqref{eq:boot-trunc} integrates over both data and weights.
\end{lemma}
\begin{proof}
\emph{Truncation and mean covariance.}
With $K=23$, Lemma~\ref{lem:proxy} gives
\[
 \|R_x\|_{L^2}\le C\sqrt n\,\alpha_m^{12}
 \le C n^{1/2-4\gamma}\le C n^{-3/2}\qquad(\|x\|\le1).
\]
Hold the suffix data fixed. Iterating the affine recursion with
$\E w_j=1$ gives $\E[e_k^{m,x,b}\mid\mathcal Z_+]=e_k^{m,x}$.
In particular $\E[F_x^b\mid\mathcal Z_+]=F_x^0$. Total covariance gives
\[
 \E\Cov(F_x^b\mid\mathcal Z_+)=\Cov(F_x^b)-\Cov(F_x^0).
\]
Under the joint law, Lemma~\ref{lem:suffixcov} bounds
$\|\Cov(F_x^b)-2\Sigma_n\|+\|\Cov(F_x^0)-\Sigma_n\|$ by
$C\{m/n+n^{-\gamma}+(n\alpha_m)^{-1}\}=O(n^{-1/2})$ for $\|x\|\le1$.
To replace $F_x^b$ by $P_x$, write
$\widetilde Y=Y-\E[Y\mid\mathcal Z_+]$. Then
$\|\widetilde R_x\|_{L^2}\le\|R_x\|_{L^2}$ and
$\|\widetilde F_x^b\|_{L^2}^2\le\tr\Cov(F_x^b)\le C$. Cauchy--Schwarz gives
\[
 \|\E V_x-\E\Cov(F_x^b\mid\mathcal Z_+)\|
 \le 2\|\widetilde F_x^b\|_{L^2}\|\widetilde R_x\|_{L^2}
        +\|\widetilde R_x\|_{L^2}^2 .
\]
Hence $\sup_{\|x\|\le1}\|\E V_x-\Sigma_n\|\le Cn^{-1/2}$.

\emph{Fluctuations with respect to the data.}
At fixed suffix data, apply Lemma~\ref{lem:covid} to $P_x$ as a
function of the $N$ independent weights at times $j\in I_+$.
Choose $I_*$ uniformly from those times and let $\mathcal I$ contain the
times preceding $I_*$ in an independent uniform permutation of $I_+$.
Use the differences in \eqref{eq:resample-differences}, with independent
copies of the weights and superscript $w$ to indicate weight replacement.
The conditional covariance identity gives
\[
 V_x=\frac N2\E\!\left[
  \Delta_{I_*}^w P_x(w^{\mathcal I})\,
  \Delta_{I_*}^w P_x(w)^\top\,\middle|\,\mathcal Z_+\right].
\]
Now change just observation $Z_i$, keeping all the weights and their
copies fixed. The first-difference bound of Lemma~\ref{lem:proxy} is
$C/\sqrt n$. If $i\ne j$, changing $Z_i$ changes a $j$th weight
difference by at most $Ch_{ij}/\sqrt n$, by the mixed-difference bound
on two distinct weighted input pairs. For $i=j$ the triangle inequality
gives $2C/\sqrt n$, which can be included by taking $h_{ii}=1$.
For each product entry and $I_*=j$ the change is therefore at most
$C h_{ij}/n$. Averaging over $I_*$ yields, uniformly in $x$ in the unit ball,
\[
 \left|\Delta_i^Z (V_x)_{ab}\right|
 \le \frac N2\,\frac1N\sum_{j\in I_+}\frac{C h_{ij}}{n}
 \le C/n .
\]
Lemma~\ref{lem:resamplingvariance} now bounds the variance of each entry by
$C N/n^2\le C/n$ and its tail about its expectation by
$2\exp(-cnt^2)$ at deviation $t$.

\emph{Uniformity in the restart.}
By Lemma~\ref{lem:proxy}, $P_x$ is affine in $x$, so
each entry of $V_x-\E V_x$ is a quadratic polynomial in $x$.
If $v_i$ are the standard basis vectors, its coefficients are determined
by the values at $0$, $\pm v_i/2$ and $(v_i+v_j)/2$ for $i<j$.
Thus its supremum on the unit ball is at most a fixed multiple,
depending on $d$, of the largest absolute value at these points.
Applying the entrywise bounds at the finitely many points followed by
the deterministic bias bound gives \eqref{eq:bootcov-rev}.
\end{proof}

\paragraph{2. The two conditioning steps.}
Return to the actual restart $x=e_m^b$ and put
$\mathcal H=\sigma(\mathcal Z_n,w_1,\ldots,w_m)$.
Conditional on $\mathcal H$, the data and initial weights are fixed.
The remaining weights are independent. Define $M=\E[U_n^b\mid\mathcal H]$.
Iterating the weighted recursion with the entire data fixed gives
$\E^b e_k^b=e_k$. It follows that
\begin{align*}
 \E^b M&=0,\\
 M&=\frac1{\sqrt n}\sum_{k=0}^{m-1}(e_k^b-e_k)
     +\frac1{\sqrt n}\sum_{k=m}^{n-1}\Gamma_{m+1:k}(e_m^b-e_m),\\
 U_n^b-M&=F_x^b-\E[F_x^b\mid\mathcal H]\qquad(x=e_m^b).
\end{align*}
In the second line the conditional mean of the weighted suffix follows the
unweighted recursion from $e_m^b$, whereas the original suffix starts at $e_m$.

For fixed even $q\ge2$, Lemma~\ref{lem:moments} bounds the $L^q$ norm
of the first sum (including its normalization) by $C_q\sqrt{m/n}$.
For the second sum, the matrix product depends only on the suffix data
and is independent of $e_m^b-e_m$ under the joint law. The product
estimate and Lemma~\ref{lem:kernel} give
\[
 \left\|\sum_{k=m}^{n-1}\Gamma_{m+1:k}\right\|_{L^q}
       \le C_q/\alpha_m,\qquad
 \|e_m^b-e_m\|_{L^q}\le C_q\sqrt{\alpha_m}.
\]
Since $\alpha_m^{-1}\lesssim m^\gamma\le m$ for $m\ge1$, we obtain
\begin{equation}\label{eq:bootshift}
 \|M\|_{L^q}\le C_q\bigl(\sqrt{m/n}+(n\alpha_m)^{-1/2}\bigr)
 \le C_q\sqrt{m/n}.
\end{equation}

\paragraph{3. The conditional normal approximation.}
Let $\lambda_n=\lambda_{\min}(\Sigma_n)$. On the data event
$H_n\le\lambda_n/2$, every $V_x$ with $\|x\|\le1$ satisfies
$\lambda_{\min}(V_x)\ge\lambda_n/2$. Given $\mathcal H$ on
$\{\|x\|\le1\}$, center $P_x$ by $\E[P_x\mid\mathcal H]$.
The covariance of this centered variable is $V_x$ evaluated at $x=e_m^b$.
There are exactly $N=n-1-m$ independent inputs in this application of
Lemma~\ref{lem:sensitivity}. Lemma~\ref{lem:proxy} bounds first and
mixed differences by $Cn^{-1/2}$ and $Ch_{ij}n^{-1/2}$ because the
original sum is normalized by $\sqrt n$. These satisfy the lemma's
$N^{-1/2}$ bounds since $N\le n$, and $N^{-1/2}=O(n^{-1/2})$ since
$N/n\to1$. The lemma, followed by Gaussian covariance comparison
(Lemma~\ref{lem:gausstools}), yields
\[
 \dc\!\left(\law(P_x-\E[P_x\mid\mathcal H]\mid\mathcal H),
                   \N(0,\Sigma_n)\right)
 \le C(n^{-1/2}+H_n).
\]

For the transfer to the exact statistic set $r=n^{-1/2}$,
$E_m=\{\|e_m^b\|\le1\}$, and define the $\mathcal Z_n$-measurable
probabilities
\[
 p_1=\Pp^b(E_m^c),\qquad
 p_2=\Pp^b\{E_m,\ \|R_x-\E[R_x\mid\mathcal H]\|>r\}.
\]
Using Lemma~\ref{lem:moments} with order $36$, then
\eqref{eq:boot-trunc} and Markov's inequality, gives
\begin{align*}
 \E p_1&\le \E\|e_m^b\|^{36}\le C\alpha_m^{18}=O(n^{-3}),\\
 \E p_2&\le r^{-2}\E[\mathbf1_{E_m}\|R_x-\E[R_x\mid\mathcal H]\|^2]
       \le r^{-2}\E[\mathbf1_{E_m}\|R_x\|^2]\le Cn^{-2}.
\end{align*}
In the last bound, condition on the data and weights through $m$.
The restart $x$ is then fixed, and the joint suffix law is precisely
the one in \eqref{eq:boot-trunc}.

Fix a deterministic convex set $B$ and set $g_B(y)=\Pp(\Sigma_n^{1/2}\eta+y\in B)$,
with $\eta$ independent of the data and weights. Conditional on $\mathcal H$,
the translation $M$ is fixed. Apply the comparison to the enlarged and eroded
sets $(B-M)^r$ and $(B-M)^{-r}$ of Lemma~\ref{lem:gausstools}. With
$\widetilde R_x=R_x-\E[R_x\mid\mathcal H]$, the bound on the Gaussian
boundary strip yields, on $\{H_n\le\lambda_n/2\}\cap E_m$,
\[
 \left|\Pp(U_n^b\in B\mid\mathcal H)-g_B(M)\right|
 \le C(n^{-1/2}+H_n+r)+\Pp(\|\widetilde R_x\|>r\mid\mathcal H).
\]
On $E_m^c$ the absolute difference is at most one. Taking conditional
expectations given $\mathcal Z_n$, using the definitions of $p_1,p_2$ and
$r=n^{-1/2}$, gives on $\{H_n\le\lambda_n/2\}$
\[
 \left|\Pp^b(U_n^b\in B)-\E^b g_B(M)\right|
       \le C(n^{-1/2}+H_n+p_1+p_2).
\]
By $\E^bM=0$ and the
uniform Hessian bound of Lemma~\ref{lem:gausstools},
\[
  |\E^b g_B(M)-g_B(0)|
    \le \tfrac12\sup_y\|\nabla^2 g_B(y)\|\,\E^b\|M\|^2
    \le C\E^b\|M\|^2 .
\]
Combine with Theorem~\ref{thm:finite} and take the supremum over $B$ to get
\begin{equation}\label{eq:bootstrap-aggregate}
 D_n\le C\big(n^{-1/2}+H_n+p_1+p_2+\E^b\|M\|^2\big).
\end{equation}
On $H_n>\lambda_n/2$ the bound follows by increasing $C$, since
$D_n\le1$ and $\inf_{n\text{ large}}\lambda_n>0$.

\paragraph{4. Mean and high-probability bounds.}
Taking expectations in \eqref{eq:bootstrap-aggregate} and using
\eqref{eq:bootcov-rev} and \eqref{eq:bootshift} with $q=2$ gives
$\E D_n=O(n^{-1/2})$. For the second assertion, Markov's inequality
and conditional Jensen give
\begin{align*}
 \Pp(p_1+p_2>n^{-1/2})&\le\sqrt n\,\E(p_1+p_2)=O(n^{-3/2}),\\
 \Pp\{\E^b\|M\|^2> n^{-1/2}\}
   &\le n^6\,\E[(\E^b\|M\|^2)^{12}]
    \le n^6\,\E\|M\|^{24}
    \le Cn^6(m/n)^{12}\le Cn^{-2}.
\end{align*}
Together with \eqref{eq:bootcov-rev}, the total exceptional probability
is $O(n^{-3/2})<1/n$ for large $n$. Inserting these bounds in
\eqref{eq:bootstrap-aggregate} proves the second assertion.\hfill$\square$

\subsection{Proof of Corollary~\ref{cor:bootcoverage}}
\begin{proof}
\emph{Bracket the bootstrap quantile.}
Write $p=1-\delta$, $s_v^2=v^\top\Sigma_n v$ and let $Z_v\sim\N(0,s_v^2)$.
Theorem~\ref{thm:finite} applied to halfspaces gives
$d_{\mathrm K}(v^\top T_n,Z_v)\le Cn^{-1/2}$.
On the data event of Theorem~\ref{thm:boot}, which has probability
at least $1-1/n$, the triangle inequality gives
\[
 d_{\mathrm K}\big(\law(v^\top U_n^b\mid\mathcal Z_n),\law(Z_v)\big)
   \le D_n+d_{\mathrm K}(v^\top T_n,Z_v)
   \le a_n:=C\sqrt{\log n/n}.
\]
For $n$ large enough that $0<p-a_n<p+a_n<1$, set
$t_\pm=s_v\Phi^{-1}(p\pm a_n)$. At every $t<t_-$ the bootstrap
distribution function is strictly below $p$, and at $t_+$ it is at least
$p$. Hence $t_-\le q_p^b(v)\le t_+$ on this event.

\smallskip\noindent\emph{Evaluate coverage at deterministic thresholds.}
The bracket gives
\[
 \Pp(v^\top T_n\le t_-)-1/n
 \le \Pp(v^\top T_n\le q_p^b(v))
 \le \Pp(v^\top T_n\le t_+)+1/n.
\]
The Gaussian probabilities at $t_\pm$ are $p\pm a_n$.
Applying the first Kolmogorov bound at each threshold gives the claimed
coverage error.
\end{proof}

\end{document}